\documentclass[11pt]{amsart}
\usepackage{amsmath,amssymb,amsfonts,amscd,graphicx}

\usepackage{latexsym}
\usepackage{amsmath}
\usepackage{xypic}
\usepackage[all]{xy}

\usepackage{dbnsymb}

\CompileMatrices

\newtheorem{theorem}{Theorem}

\numberwithin{equation}{section}

\newtheorem{cor}{Corollary}[section]

\newtheorem{lem}[cor]{Lemma}
\theoremstyle{definition}
\newtheorem{defi}[cor]{Definition}

\newcommand{\co}{\colon\thinspace}

\begin{document}

\title[Tutte chromatic identities from the Temperley-Lieb algebra]
{Tutte chromatic identities\\ from the Temperley-Lieb algebra}
%\shorttitle{Tutte chromatic identities}

\author[Paul Fendley and Vyacheslav Krushkal]{Paul Fendley$^1$ and Vyacheslav Krushkal$^2$}
\address{Paul Fendley, Department of Physics, University of Virginia,
Charlottesville, VA 22904 USA;
and All Souls College and the Rudolf Peierls Centre for Theoretical Physics,
University of Oxford, 1 Keble Road, Oxford OX13NP, UK
}

\email{fendley\char 64 virginia.edu}

\address{Vyacheslav Krushkal, Department of Mathematics, University of Virginia,
Charlottesville, VA 22904-4137 USA}

\email{krushkal\char 64 virginia.edu}

\thanks{$^1$ Supported in part by NSF grants DMR-0412956 and DMR/MSPA-0704666, and by the UK EPSRC under grant EP/F008880/1}
\thanks{$^2$ Supported in part by NSF grants DMS-0605280 and DMS-0729032, and by Harvard University}
%\thanks{Date: October 24, 2007}

\begin{abstract}
This paper introduces a conceptual framework, in the context of quantum topology
and the algebras underlying it, for analyzing relations
obeyed by the chromatic polynomial $\chi(Q)$ of planar graphs.
Using it we give new proofs and substantially extend a number of classical
results concerning the combinatorics of the chromatic polynomial.
In particular, we show that
Tutte's golden identity is a consequence of level-rank
duality for $SO(N)$ topological quantum field theories and
Birman-Murakami-Wenzl algebras. This identity is a remarkable feature of the chromatic
polynomial relating
$\chi({\phi+2})$ for any triangulation of the sphere to
$(\chi({\phi+1}))^2$ for the same graph, where ${\phi}$ denotes the
golden ratio. The new viewpoint presented here explains that Tutte's identity
is special to these values of the parameter $Q$.
A natural context for analyzing such properties of the
chromatic polynomial is provided by the {\em chromatic algebra},
whose Markov trace is the chromatic polynomial of an associated
graph.  We use it to show that another identity of
Tutte's for the chromatic polynomial at $Q={\phi}+1$ arises from a
Jones-Wenzl projector in the Temperley-Lieb algebra. We generalize
this identity to each value $Q= 2+2\cos(2\pi j/(n+1))$ for $j< n$
positive integers.  When $j=1$, these $Q$ are the Beraha numbers,
where the existence of such identities was conjectured by Tutte.  We
present a recursive formula for this sequence of chromatic polynomial
relations.
\end{abstract}

\maketitle

\section{Introduction}

In a series of papers in 1969 (cf \cite{T1}, \cite{T2}), W.T. Tutte
discovered several remarkable properties of the chromatic polynomial $\chi(Q)$
of planar graphs evaluated at the special value $Q={\phi}+1$,
where ${\phi}$ denotes the golden ratio, ${\phi}=\frac{1+\sqrt5}{2}$.
Among these is the ``golden identity'': for a planar triangulation $T$,
\begin{equation} \label{Tutte's identity1}
  {\chi}^{}_T({\phi}+2)=({\phi}+2)\; {\phi}^{3\,V(T)-10}\,
  ({\chi}^{}_T({\phi}+1))^2,
\end{equation}
where $V(T)$ is the number of vertices of the triangulation. Tutte used this identity
to establish that ${\chi}^{}_T({\phi}+2)$ is positive, a result interesting in connection
to the four-color theorem. Another
property (see (11.15) in \cite{T3}) is the relation
\begin{equation}
\label{Tutte's identity3}
{\chi}^{}_{Z_1}({\phi}+1)+{\chi}^{}_{Z_2}({\phi}+1)=
{\phi}^{-3}[{\chi}^{}_{Y_1}(\phi+1)+{\chi}^{}_{Y_2}({\phi}+1)],\end{equation}
where $Y_i$, $Z_i$ are planar graphs which are
locally related as shown in figure 1.

The main purpose of this paper is to show that Tutte's results
naturally fit in the framework of quantum topology and the algebras
underlying it. We give conceptual proofs of these identities, very
different from the combinatorial proofs of Tutte's. This allows us to
show that both identities (\ref{Tutte's identity1},\ref{Tutte's identity3})
are not {\it ad hoc}, but are (particularly elegant) consequences of a
deeper structure. We generalize the linear relation (\ref{Tutte's
identity3}) to other special values of $Q$, dense in the interval
$[0,4]$.

For example, the identity analogous to (\ref{Tutte's identity3}) at
$Q=2$ is simply that $\chi_\Gamma^{}(2)=0$ whenever the graph $\Gamma$
includes a triangle. This vanishing is completely obvious from the
original definition of $\chi_\Gamma(Q)$ for integer $Q$ as the number
of $Q$-colorings of $\Gamma$ such that adjacent vertices are colored
differently. Nevertheless, this vanishing does not follow immediately
from the contraction-deletion relation used to define $\chi(Q)$ for
all $Q$. We use the algebraic approach to fit this obvious identity
and the rather non-trivial identity (\ref{Tutte's identity3}) into a
series of identities at special values of $Q$ which are independent of
the contraction-deletion relation, although of course are consistent
with it.

We show that Tutte's golden identity (\ref{Tutte's identity1}) is a
consequence of level-rank duality for $SO(N)$ TQFTs and the
Birman-Murakami-Wenzl algebras underlying the construction of the
(doubled) TQFTs.  Level-rank duality is an important property of
conformal field theories and topological quantum field theories. It
implies that the $SO(N)$ level $k$, and the $SO(k)$ level $N$ theories
are isomorphic, cf.\ \cite{MNRS} (in our conventions the level $k$ of
$SO(3)$ corresponds to the level $2k$ of $SU(2)$, so the $SO(3)$
theories are labeled by the half-integers). We exploit the isomorphism
between the $SO(3)_4$ and $SO(4)_3$ theories, and then show that the
latter splits into a product of two copies of $SO(3)_{3/2}$. The
partition function of an $SO(3)$ theory is given in terms of the
chromatic polynomial, specifically ${\chi}({\phi}+2)$ for $SO(3)_4$
and ${\chi}({\phi}+1)$ for $SO(3)_{3/2}$. These
isomorphisms enable us to give a rigorous proof
of the golden identity (\ref{Tutte's identity1}). This new viewpoint
makes it clear that the fact that golden identity relates a chromatic
polynomial squared to the chromatic polynomial at another value is
very special to these values of $Q$: generalizations using level-rank
duality $SO(N)_4\leftrightarrow SO(4)_N$ do not involve the chromatic
polynomial.

An important and rather natural tool for analyzing the chromatic
polynomial relations is the {\em chromatic algebra} ${\mathcal C}^Q_n$
introduced (in a different context) in \cite{Martin}, and analyzed in
depth in \cite{FK}.  The basic idea in the definition of the chromatic
algebra is to consider the contraction-deletion rule as a linear
relation in the vector space spanned by graphs, rather than just a
relation defining the chromatic polynomial. In this context, the right
object to study is the space of {\em dual} graphs, and the algebra is
defined so that the Markov trace of a graph is the chromatic
polynomial of its dual.  The parameter $Q$ in the chromatic algebra
${\mathcal C}^Q$ is related to the value of the loop (or dually the
value of the chromatic polynomial of a single point.)

Relations between algebras and topological invariants like the
chromatic polynomial are familiar from both knot theory and
statistical mechanics. Our companion paper \cite{FK} discusses
these connections in detail. In particular, the relation
between the $SO(3)$ Birman-Wenzl-Murakami algebra and the chromatic
algebra described initially in \cite{FR} is derived there. This in
turn yields a relation between the chromatic algebra and the
Temperley-Lieb algebra \cite{FR}, which is described in section
\ref{sec:TL} below. Readers familiar with the Potts model and/or the Tutte
polynomial may recall that it has long been known (indeed from
Temperley and Lieb's original paper \cite{TL}) that the chromatic
polynomial and the more general Tutte polynomial can be related to the
Markov trace of elements of the Temperley-Lieb algebra. The relation
derived in \cite{FR} and here is quite different from the earlier
relation; in statistical-mechanics language ours arises from the
low-temperature expansion of the Potts model, while the earlier one
arises from Fortuin-Kasteleyn cluster expansion \cite{FK}.

The relation of the chromatic algebra to the Temperley-Lieb algebra
utilized here makes it possible to rederive and greatly extend Tutte's
identities. In this context, linear identities such as Tutte's
(\ref{Tutte's identity3}) are understood as finding elements of the
{\em trace radical}: elements of the chromatic algebra which,
multiplied by any other element of the algebra, are in the kernel of
the Markov trace.  For example, it is well known how the Jones-Wenzl
projectors $P^{(k)}$, crucial for the construction of $SU(2)$
topological quantum field theories, are defined within the
Temperley-Lieb algebra. We define these relations in the chromatic
algebra, for special values of the parameter $Q$, as a pull-back of
the Jones-Wenzl projectors.  Tutte's relation (\ref{Tutte's
identity3}) is a pull-back of the projector $P^{(4)}$, corresponding
to $Q=\phi+1$.  More precisely, in this algebraic setting the relation
(\ref{Tutte's identity3}) for the chromatic polynomials of the graphs
in figure 1 is a consequence of the fact that
\begin{equation}
\label{Tutte's identity4}\widehat Z_1\, +\, \widehat Z_2\; -
{\phi}^{-3}\,[\,\widehat Y_1\, +\, \widehat Y_2\,]
\end{equation}
is an element of the trace radical of the algebra ${\mathcal C}^{{\phi}+1}_2$.
Here $\widehat Y_i$, $\widehat Z_i$ are the planar graphs dual to the
graphs $Y_i, Z_i$, see figure \ref{fig:dualTutte}.

\begin{figure}[t]
%\centering
\includegraphics[height=2.2cm]{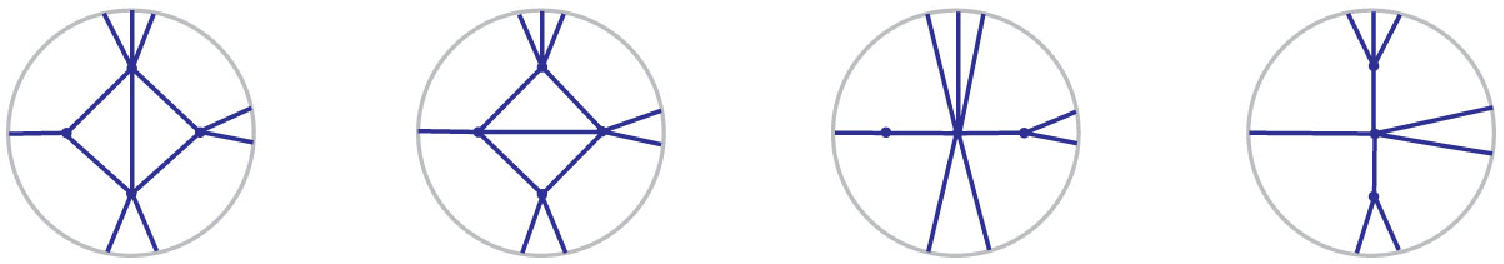}
{\scriptsize
    \put(-365,3){$Z_1$}
    \put(-265,3){$Z_2$}
    \put(-166,3){$Y_1$}
    \put(-65,3){$Y_2$}}
\caption{The graphs in Tutte's identity (\ref{Tutte's identity3}). There
  can be any number of lines from each vertex to the boundary of the
  disk, and the graphs are identical outside the disk.}
\label{fig:Tutte3}
\end{figure}

\begin{figure}[ht]
%\centering
\includegraphics[height=2.2cm]{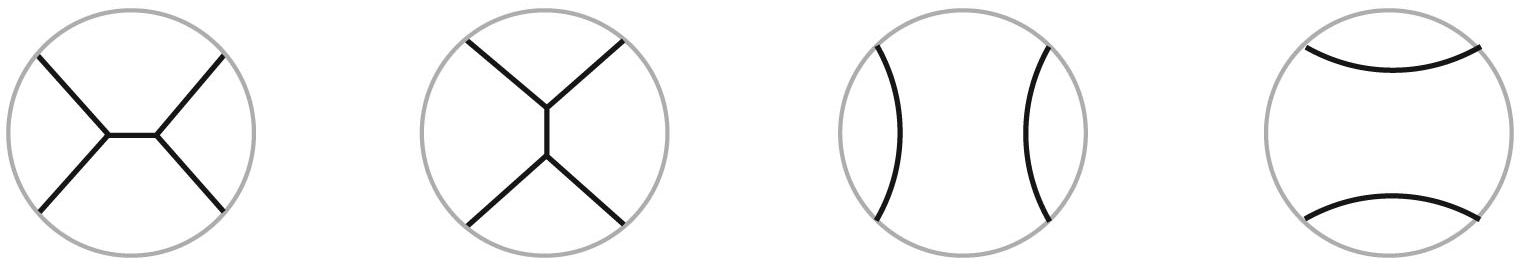}
{\scriptsize
    \put(-369,3){$\widehat Z_1$}
    \put(-269,3){$\widehat Z_2$}
    \put(-169,3){$\widehat Y_1$}
    \put(-68,3){$\widehat Y_2$}}
\caption{The graphs $\widehat Y_i$, $\widehat Z_i$ dual to the graphs
in figure 1.}
\label{fig:dualTutte}
\end{figure}

Tutte remarks about (\ref{Tutte's identity3}): ``this equation can be
taken as the basic one in the theory of golden chromials''
\cite{T3}. This statement has a precise meaning in our algebraic
context: $P^{(4)}$ is a generator of the unique proper ideal of the
Temperley-Lieb category related to the chromatic polynomial at
$Q=\phi+1$. Note that Tutte's discovery of this relation (for this specific value of $Q$)
predates that
of Jones and Wenzl by fifteen years!

We show that Tutte's identity (\ref{Tutte's identity3}) has an analog
 for any value of $Q$ obeying
\begin{equation}
Q=2+2\cos\left(\frac{2\pi j}{n+1}\right)
\label{QJW}
\end{equation}
for $j$ and $n$ positive integers obeying $j< n$.  The requirement of
integer $n$ arises from the structure of the Temperley-Lieb category:
these are the only values for which the Temperley-Lieb category has a
non-trivial proper ideal (given by the trace radical). When $j=1$,
these values of $Q$ are the Beraha numbers $B_{n+1}$, discussed in
more detail below. Tutte conjectured in \cite{T3} that there is such a
relation, similar to (\ref{Tutte's identity3}) at ${\phi}+1=B_5$, for
each Beraha number, and we give a recursive formula for these
relations based on the formula for the Jones-Wenzl projectors
$P^{(n)}$.  In fact, based on the relation between the chromatic
algebra and the $SO(3)$ BMW algebra \cite{FK}, it seems reasonable to
conjecture that these are {\em all} linear relations which preserve
the chromatic polynomial at a given value of $Q$ (and in particular
there are no such relations at $Q$ not equal to one of special values
(\ref{QJW}).)

We use similar ideas to give also a direct algebraic proof of the golden
identity (\ref{Tutte's identity1}).
This proof (given in section \ref{golden section}) is motivated by the argument using
level-rank duality, discussed above,
but it is given entirely in the context of the chromatic algebra.
It is based on a map ${\Psi}\co
{\mathcal C}^{{\phi}+2}\longrightarrow ({\mathcal
C}^{{\phi}+1}/R)\times ({\mathcal C}^{{\phi}+1}/R)$, where $R$ denotes
the trace radical.  We show that the existence of this map is implied
by the relation (\ref{Tutte's identity3}), and then (\ref{Tutte's
identity1}) follows from applying the algebra traces to the
homomorphism ${\Psi}$.
%This contrasts with the linear relation (\ref{Tutte's identity3})
%which we generalize to other special values of $Q$.
The golden identity has an interesting application in physics in
quantum loop models of ``Fibonacci anyons'', where it implies that
these loop models should yield topological quantum field theories in
the continuum limit (see \cite{Fnew,FFNWW,Walker} for more details.)

We conclude this general discussion of Tutte's results by noting his estimate:
\begin{equation} \label{Tutte's identity5} |{\chi}^{}_T({\phi}+1)|\leq
{\phi}^{5-k},\end{equation} where $T$ is a planar triangulation and
$k$ is the number of its vertices. The origin of Beraha's definition
of the numbers $B_n$ is his observation \cite{Be} that the zeros of
the chromatic polynomial of large planar triangulations seem to
accumulate near them. Tutte's estimate (\ref{Tutte's identity5})
gives a hint about this phenomenon for $B_5$. Efforts have been made (see
\cite{KS}, \cite{S}) to explain Beraha's observation using quantum
groups. We hope that our approach will shed new light on this
question.

The paper is organized as follows.  Section \ref{sec:TL} introduces the
chromatic algebra ${\mathcal C}_n^Q$, reviews the standard material on
the Temperley-Lieb algebra $TL_n^d$, and describes the algebra
homomorphism ${\mathcal C}_n^Q\longrightarrow TL_{2n}^{\sqrt Q}$.  It
then shows that Tutte's relation (\ref{Tutte's identity3}) is the
pullback of the Jones-Wenzl projector $P^{(4)}\in TL^{\phi}$. With this
relation at hand, section \ref{golden section} gives a direct
algebraic proof of Tutte's golden identity (\ref{Tutte's identity1}) in the chromatic
algebra setting. Section \ref{sec:levelrank} discusses the level-rank duality
of $SO(N)$ BMW algebras, and shows that the golden identity (\ref{Tutte's identity1})
is its consequence.
In the final section \ref{Beraha} we show that for any $Q$ obeying
(\ref{QJW}) the chromatic polynomial obeys a generalization of
(\ref{Tutte's identity3}). We give a recursive formula for this
sequence of chromatic relations.  This paper is largely
self-contained, but we refer the reader interested in the
chromatic algebra, the $SO(3)$ Birman-Murakami-Wenzl algebra, and the
relations between the chromatic polynomial, link invariants and TQFTs,
to our companion paper \cite{FK}.

%\medskip

\section{The chromatic algebra and the Temperley-Lieb algebra} \label{sec:TL}

The definition of the algebraic structure of the chromatic algebra
${\mathcal C}^Q_n$ -- the product structure, the trace, the inner
product -- is motivated by that of the Temperley-Lieb algebra.
In this section we set up the general framework for relating the chromatic
polynomial to the Temperley-Lieb algebra, and we show that the relation
(\ref{Tutte's identity3}) corresponds to a Jones-Wenzl projector.
We start by defining the chromatic algebra \cite{Martin,FK}, and then
review standard material on Temperley-Lieb
algebra (the reader is referred to \cite{KL} for more details),
We then define an algebra homomorphism ${\mathcal C}^Q_n\longrightarrow
TL^d_{2n}$, where $Q=d^2$, which respects these structures. In
particular, the pullback under this homomorphism of the {\em trace
radical} in $TL^d_{2n}$ (the ideal consisting of elements whose inner
product with all other elements in the algebra is trivial) is in the
trace radical of ${\mathcal C}^Q_{2n}$, which corresponds to local
relations on graphs which preserve the chromatic polynomial of the
duals, for a given value of $Q$.
In fact, the relevant algebraic
structure here is the chromatic, respectively Temperley-Lieb, {\em
category}.

\subsection{The chromatic algebra}
The {\em chromatic polynomial} ${\chi}^{}_{\Gamma}(Q)$ of a graph
$\Gamma$, for $Q\in {\mathbb Z}_+$, is the number of colorings of the
vertices of $\Gamma$ with the colors $1,\ldots, Q$ where no two
adjacent vertices have the same color. A basic property of the
chromatic polynomial is the contraction-deletion rule: given any edge
$e$ of $\Gamma$ which is not a loop,
\begin{equation} \label{chromatic poly1}
{\chi}^{}_{\Gamma}(Q)={\chi}^{}_{{\Gamma}\backslash e}(Q)-{\chi}^{}_{{\Gamma}/e}(Q)
\end{equation}
where ${\Gamma}\backslash e$ is the graph obtained from $\Gamma$ by
deleting $e$, and ${\Gamma}/e$ is obtained from $\Gamma$ by
contracting $e$. If $\Gamma$ contains a loop then
${\chi}^{}_{\Gamma}\equiv 0$; if $\Gamma$ has no edges and $V$ vertices, then
$\chi^{}_\Gamma(Q)=Q^V$. These properties enable one to define
${\chi}^{}_{\Gamma}(Q)$ for any, not necessarily integer, values of $Q$.

\begin{figure}[ht]
%\centering
%\vspace{.2cm}
\includegraphics[width=2.8cm]{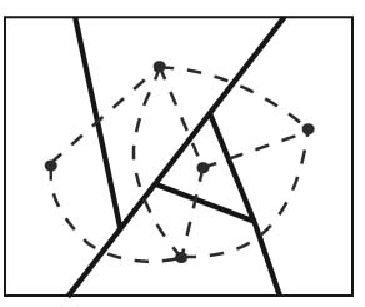}
{\scriptsize
    \put(-91,5){$R$}
    \put(-71,53){$G$}
    \put(-16,45){$\widehat G$}}
 \caption{A basis element $G$ of the algebra ${\mathcal F}_2$ and the dual graph
$\widehat G$ (drawn dashed).}
\label{fig:GGdual}
\end{figure}

Consider the free algebra ${\mathcal F}_n$ over ${\mathbb C}[Q]$
whose elements are formal
linear combinations of the isotopy classes of trivalent graphs in a
rectangle $R$ (figure \ref{fig:GGdual}). The intersection of each such graph with the boundary
of $R$ consists of precisely $2n$ points: $n$ points at the top and
the bottom each, and the isotopy, defining equivalent graphs, is
required to preserve the boundary.
Note that the vertices of the
graphs in the interior of $G$ are trivalent, in particular they do not have ends
($1$-valent vertices) other than those on the boundary of $R$. It is convenient
to allow $2$-valent vertices as well, so there may be loops disjoint from the rest
of the graph.
The multiplication in ${\mathcal F}_n$ is given by
vertical stacking, and the inclusion of algebras
${\mathcal F}_n\subset {\mathcal F}_{n+1}$ is defined on the graphs generating
${\mathcal F}_n$ as the addition of a vertical strand on the right.
Given $G\in{\mathcal F}_n$, the vertices of its
{\em dual graph} $\widehat G$ correspond to the complementary regions
$R\smallsetminus G$, and two vertices are joined by an edge in
$\widehat G$ if and only if the corresponding regions share an edge,
as illustrated in figure \ref{fig:GGdual}.

\begin{defi} \label{chromatic definition}
The chromatic algebra in degree $n$, ${\mathcal C}_n$, is the algebra
over ${\mathbb C}[Q]$ which is defined as the quotient of the free
algebra ${\mathcal F}_n$ by the ideal $I_n$ generated by the relations
shown in figure \ref{fig:trivalent}. In addition, the value of a
trivial simple closed curve is set to be $Q-1$. When the parameter $Q$ is specialized to
a complex number, the resulting algebra over ${\mathbb C}$ is denoted
${\mathcal C}_n^Q$.
Set ${\mathcal C}=\cup_n {\mathcal C}_n$.

\begin{figure}[ht]
%\centering
%\vspace{.2cm}
\includegraphics[height=1.85cm]{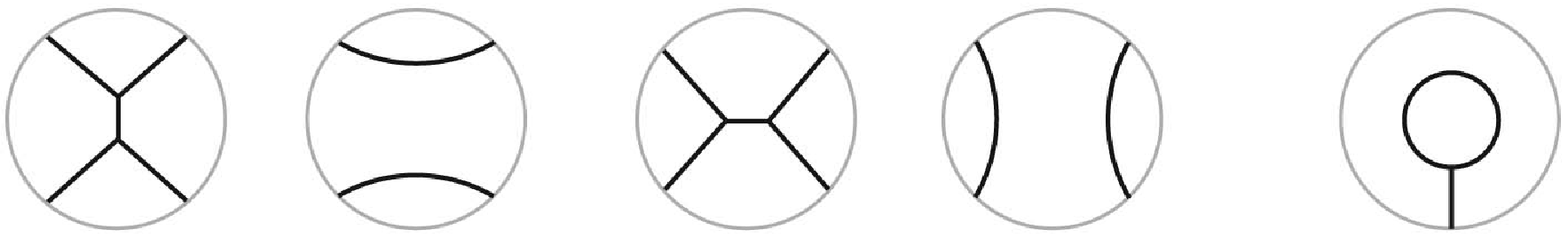}
    \put(-296,25){$+$}
    \put(-226,25){$=$}
    \put(-155,25){$+$}
    \put(4,25){$=\;0.$}
    \put(-87,8){$,$}
\caption{Relations in the trivalent presentation of the chromatic algebra.}
\label{fig:trivalent}
\end{figure}
\end{defi}
The first relation in figure \ref{fig:trivalent} is sometimes known as the ``H-I
relation'', while the second is requiring that ``tadpoles'' vanish. A few words may be
helpful in explaining the relations defining the chromatic algebra. If there is a simple
closed curve bounding a disk in a rectangular picture, disjoint from the rest of the graph, it may be erased
while the element represented in ${\mathcal C}_n$ is multiplied by $Q-1$. The ideal corresponding to the
relation on the left in figure \ref{fig:trivalent} is generated by linear combinations of graphs in ${\mathcal F}_n$
which are identical outside a disk in the rectangle, and which differ according to this relation in the disk,
figure \ref{fig:ex7}. This may be naturally expressed in the language of planar algebras (see also section
\ref{Beraha}.)

\begin{figure}[ht]
%\centering
\includegraphics[width=2.8cm]{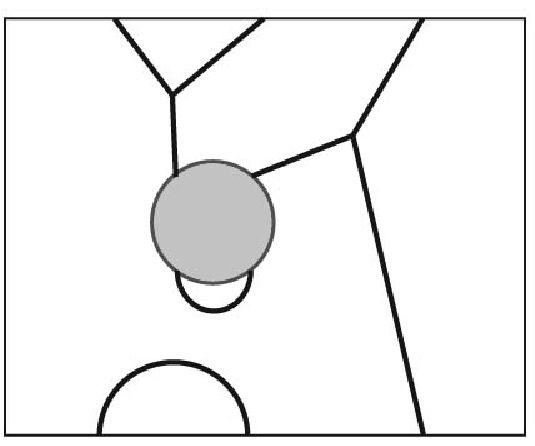} \caption{An element in the ideal $I_3$. The shaded disk contains the relation on the left
in figure \ref{fig:trivalent}.}
\label{fig:ex7}
\end{figure}

{\bf Remark.} In this definition we used just the trivalent graphs,
and this is sufficient for the purposes of this section, in particular
for the proofs of Tutte's identities (\ref{Tutte's identity1}),
(\ref{Tutte's identity3}).  The definition of the chromatic algebra
using all planar graphs (isomorphic to the trivalent one considered
here), is given in section \ref{Beraha}.

\begin{defi}
The trace, $tr_{\chi}\co {\mathcal C}^Q\longrightarrow{\mathbb C}$ is defined on the additive
generators (graphs) $G$ by connecting the
endpoints of $G$ by disjoint arcs in the complement of the rectangle $R$ in the plane (denote the result by $\overline G$)
and evaluating $$tr_{\chi}(G)\; =\; Q^{-1}\cdot {\chi}_{\widehat{\overline G}}(Q).$$
\end{defi}

\begin{figure}[h]
%\centering
\includegraphics[height=2.3cm]{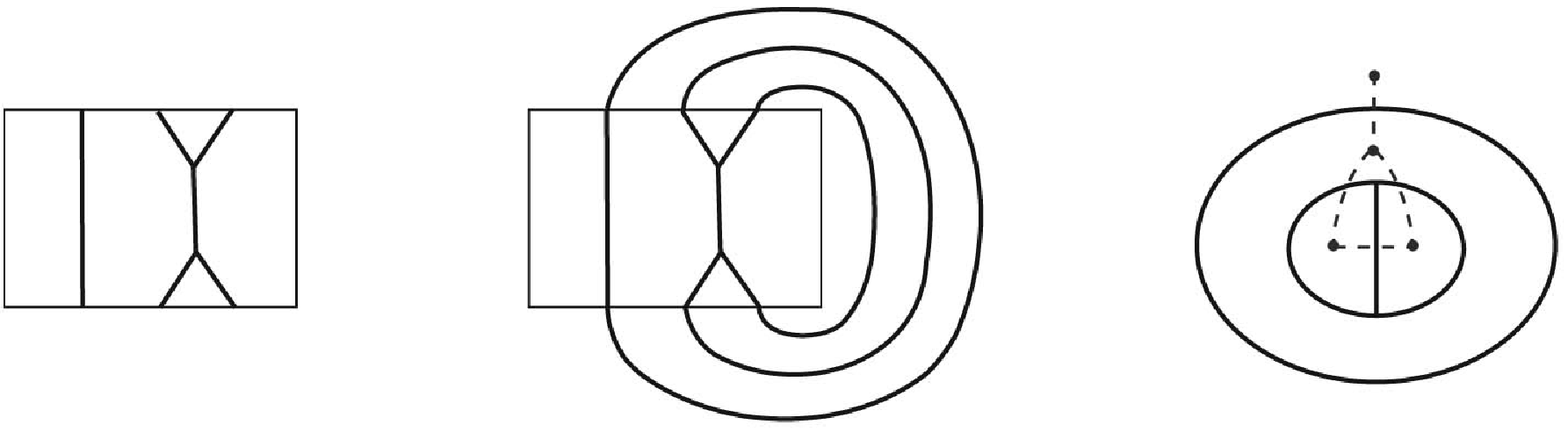}
{\small
    \put(-248,32){$Tr$}
    \put(-172,32){$=$}
    \put(-73,32){$=$}
    \put(2,32){$=(Q-1)^2(Q-2)$.}}
\caption{An example of the evaluation of the trace.}
\label{fig:ex4}
\end{figure}

Figure \ref{fig:ex4} shows the trace of an example. The factor $Q^{-1}$ provides a normalization of the trace which
turns out to be convenient from
the point of view of the relation with the Temperley-Lieb algebra (see below). One checks that the trace
is well-defined by considering the chromatic polynomial of the dual graphs. Specifically, the second
relation in figure \ref{fig:trivalent} holds since the dual graph has a loop. The first relation holds since (in the
notation in figure 1), the deletion-contraction rule for the chromatic polynomial implies
\begin{equation}
{\chi}_{Z_1}^{}(Q)+{\chi}^{}_{Y_1}(Q)={\chi}^{}_{Z_2}(Q)+{\chi}^{}_{Y_2}(Q).
\end{equation}

Finally, the relation replacing a trivial simple closed curve by a factor $(Q-1)$ corresponds to the effect
on the chromatic polynomial of the dual graph of erasing a $1$-valent vertex and of the adjacent edge.

\subsection{The Temperley-Lieb algebra} \label{TL subsection}
The Temperley-Lieb algebra in degree $n$, $TL_n$, is an algebra over
${\mathbb C}[d]$ generated by $1, E_1,\ldots, E_{n-1}$ with the
relations \cite{TL}
\begin{equation} \label{TL relations}
E_i^2=E_i, \qquad E_iE_{i\pm 1} E_i=\frac{1}{d^2}\,E_i, \qquad E_i E_j=E_j E_i \; \, {\rm for} \; \,
|i-j|>1.
\end{equation}
Define $TL=\cup_n TL_n$.
The indeterminate $d$ may be set to equal a specific complex number, and
when necessary, we will include this in the notation, $TL^d_n$.

\begin{figure}[ht]
\includegraphics[width=9.7cm]{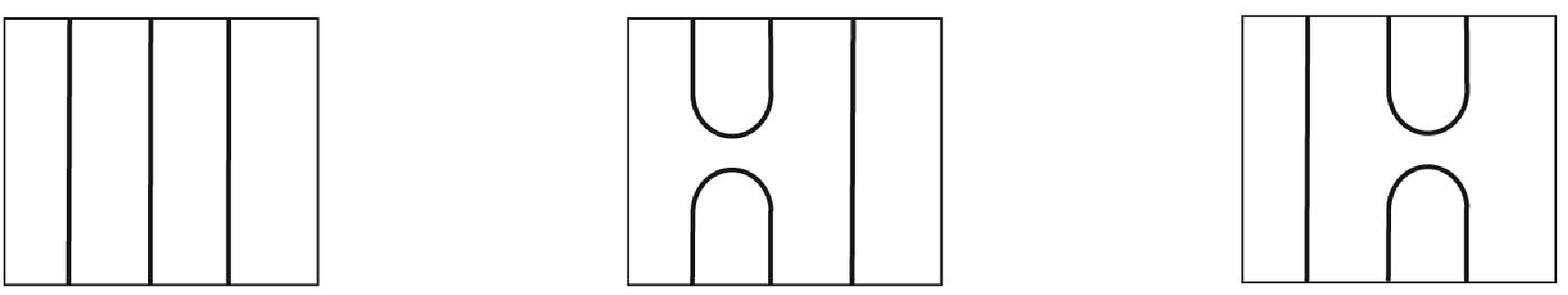}
    \put(-298,22){$1 =$}
    \put(-205,22){$E_1 =\frac{1}{d}$}
    \put(-96,22){$E_2 = \frac{1}{d}$}
\vspace{.2cm} \caption{Generators of $TL_3$}
\label{fig:TL}
\end{figure}

It is convenient to represent the elements of $TL_n$ pictorially:
in this setting, an element of $TL_n$ is a linear combination of
$1-$dimensional submanifolds in a rectangle $R$. Each submanifold meets
both the top and the bottom of the rectangle in exactly $n$ points. The multiplication then corresponds
to vertical stacking of rectangles. The generators
of $TL_3$ are illustrated in figure \ref{fig:TL}.

The rectangular pictures are considered equivalent if they are
isotopic relative to the boundary. Another equivalence arises from the
relation $E_i^2=E_i$: the element in $TL$
corresponding to a picture in $R$ with a circle (a simple closed
curve) is equivalent to the element with the circle deleted and
multiplied by $d$. Isotopy together with this relation are sometimes
referred to as {\em $d-$isotopy} \cite{F}.

The trace $tr_d\co TL^d_n\longrightarrow {\mathbb C}$ is defined on
the additive generators (rectangular pictures) by connecting the top
and bottom endpoints by disjoint arcs in the complement of $R$ in the
plane (the result is a disjoint collection of circles in the plane),
and then evaluating $d^{\# circles}$. The Hermitian product on $TL_n$
is defined by $\langle a,b\rangle=tr(a\, \bar b)$, where the
involution $\bar b$ is defined by conjugating the complex
coefficients, and on an additive generator $b$ (a picture in $R$) is
defined as the reflection in a horizontal line.

\subsection{A map from the chromatic algebra to the Temperley-Lieb algebra} \label{relations}
\begin{figure}[ht]
%\centering
\includegraphics[width=11.5cm]{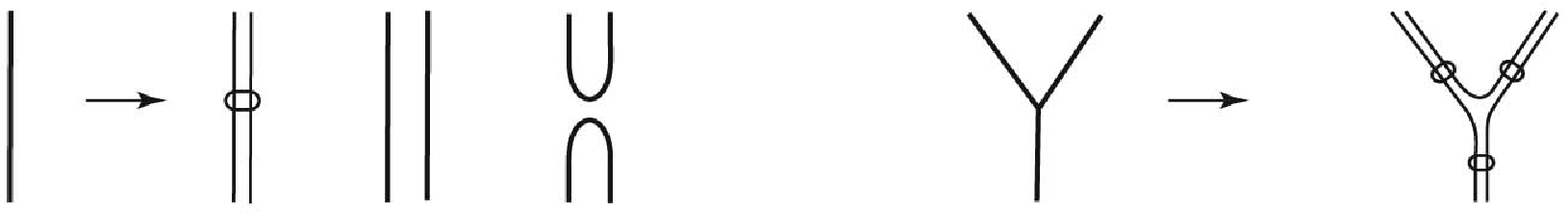}
    \put(-307,29){${\Phi}$}
    \put(-265,21){$=$}
    \put(-232,21){$-\frac{1}{d}$}
    \put(-80,29){${\Phi}$}
    \put(-57,21){$d^{1/2}\;\cdot$}
\caption{Definition of the homomorphism ${\Phi}\co {\mathcal C}^Q_n\longrightarrow TL^{\sqrt Q}_{2n}$}
\label{fig:phi1}
\end{figure}

\begin{defi}
Define a homomorphism ${\Phi}\co {\mathcal F}_n\longrightarrow TL^d_{2n}$ on the multiplicative generators
(trivalent graphs in a rectangle) of the free algebra ${\mathcal F}_n$
by replacing each edge with the linear combination ${\Phi}(\;|\;)=\smoothing-\frac{1}{d}\,\hsmoothing\,$, and resolving
each vertex as shown in figure \ref{fig:phi1}. Moreover, for a graph $G$, ${\Phi}(G)$ contains a factor
$d^{V(G)/2}$, where $V(G)$ is the number of vertices of $G$.
\end{defi}

\begin{figure}[ht]
%\centering
\includegraphics[height=1.6cm]{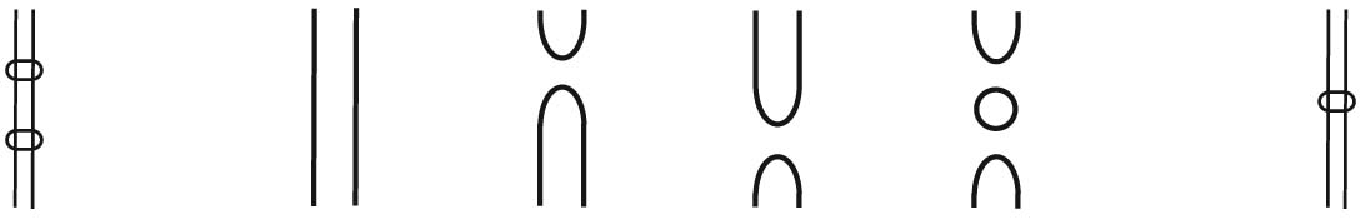}
    \put(-250,21){$=$}
    \put(-200,21){$-\;\frac{1}{d}$}
    \put(-151,21){$-\,\frac{1}{d}$}
    \put(-108,21){$+\,\frac{1}{d^2}$}
    \put(-43,21){$=$}
\caption{$P_2\circ P_2=P_2$}
\label{fig:JW}
\end{figure}

{\bf Remarks.} 1. The reader may have noticed that ${\Phi}$ replaces each edge with the
second Jones-Wenzl projector $P_2$, well-known in the study of the Temperley-Lieb algebra \cite{We}.
They are idempotents: $P_2\circ P_2=P_2$, and this identity (used below) may be easily checked directly, figure \ref{fig:JW}.

2. Various authors have considered versions of the map ${\Phi}$ in the
 knot-theoretic and TQFT contexts, see \cite{Y,J,KL, FFNWW,Walker,Fnew}. In \cite{FR}
 this was used to give a map of the SO(3) BMW algebra to
 the Temperley-Lieb algebra, probably the first instance where this
 map is considered as an algebra homomorphism.

\begin{lem} \sl ${\Phi}$ induces a well-defined homomorphism of algebras ${\mathcal C}^Q_n\longrightarrow TL^d_{2n}$,
where $Q=d^2$.
\end{lem}
To prove this lemma, one needs to check that the relations in ${\mathcal C}^Q_n$, figure \ref{fig:trivalent}, hold in the Temperley-Lieb
algebra. It follows from the definition of ${\Phi}$ (figure \ref{fig:phi1}) that

\begin{figure}[ht]
%\centering
\includegraphics[height=1.5cm]{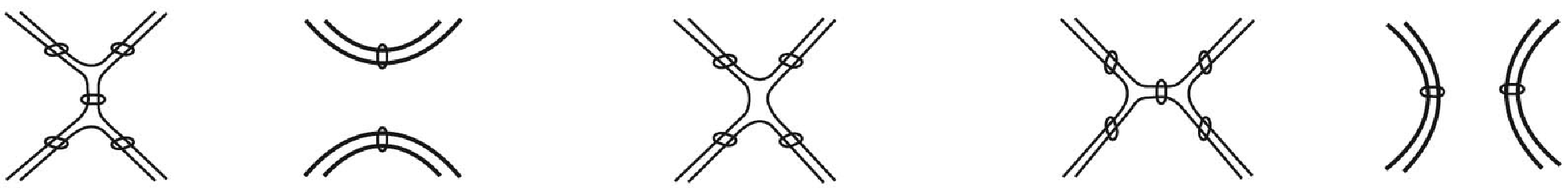}
    \put(-380,17){$d\cdot$}
    \put(-320,17){$+$}
    \put(-245,17){$=\;\;\;\; d\cdot$}
    \put(-155,17){$=\;\;\;\; d\cdot $}
    \put(-65,17){$+$}
\caption{}
\end{figure}

Similarly, one checks the other two defining relations of the chromatic algebra. \qed

The following lemma implies that the homomorphism ${\Phi}$ preserves the trace of
the chromatic, respectively Temperley-Lieb, algebras (see also Theorem 1 in \cite{FFNWW}).

\begin{lem} \label{commutative} \sl Let $G$ be a trivalent planar graph. Then
\begin{equation} \label{chromatic dual}
Q^{-1}\, {\chi}_Q(\widehat G)\, =\, {\Phi}(G).\end{equation}
 Here $Q=d^2$ and, abusing the notation, we denote by ${\Phi}(G)$ the evaluation
$d^{\#}$ applied to the linear combination of simple closed curves obtained by applying ${\Phi}$ as shown in figure
\ref{fig:phi1}.
Therefore, the following diagram
commutes:
\begin{equation} \label{chromaticTL}
\xymatrix{ {\mathcal C}_n^{\scriptsize Q}  \ar[d]^{tr_{\chi}} \ar[r]^{\Phi} & TL_{2n}^{d}  \ar[d]^{tr_{d}}\\
{\mathbb C}\ar[r]^=  &  {\mathbb C} }\end{equation}
\end{lem}

\begin{figure}[ht]
%\centering
\includegraphics[height=2.2cm]{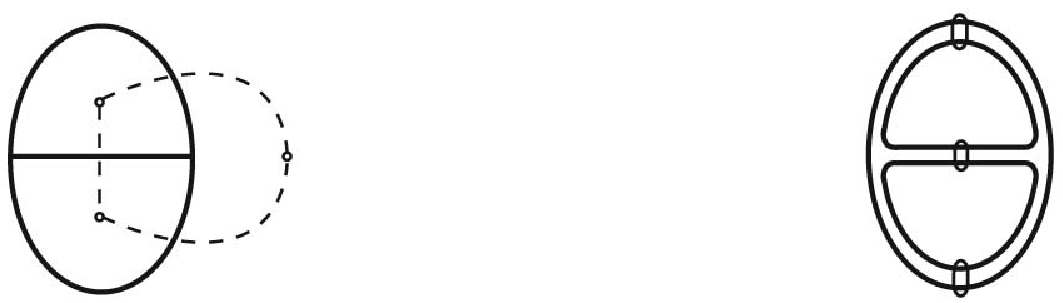}
{\small
    \put(-235,37){$G$}
    \put(-165,7){$\widehat G$}
    \put(-55,27){$d\, \cdot$}
    \put(-95,27){${\Phi}(G)\, =$}}
\caption{}
\label{fig:ex5}
\end{figure}

For example, for the theta-graph $G$ in figure \ref{fig:ex5}, one checks that $$Q^{-1}
{\chi}_Q(\widehat G)\, =\, (Q-1)(Q-2)\, =\, d^4-3d^2+2\, =\, {\Phi}(G).$$

{\em Proof of lemma \ref{commutative}.} We use the state sum formula for the chromatic polynomial (cf \cite{B}):
$${\chi}_Q(\widehat G)\, =\, \sum_{s\subseteq E(\hat
G)} (-1)^{|s|}\, Q^{k(s)} .$$ Here $k(s)$ is the number of connected
components of the graph $\widehat G_s$ whose vertex set is $V(\widehat G)$ and the
edge set is $s$.

First assume that $G$ is a connected graph. Recall that ${\Phi}(G)$ is
obtained from $G$ by replacing each edge $|$ by the linear combination
$(\smoothing - 1/d \, \hsmoothing)$, and resolving each vertex as
shown in figure \ref{fig:phi1}. Then ${\Phi}(G)$ is (the evaluation
of) a linear combination of simple closed curves. This linear
combination can be represented as a sum parametrized by the subsets of
the set of the edges of the dual graph $\hat G$, $s\subset E(\hat
G)$. For each such subset, the corresponding term is obtained by
replacing each edge $|$ of $G$ not intersecting $s$ with $\smoothing$,
and each edge intersecting $s$ with $\hsmoothing$.  Moreover, this
collection of simple closed curves is the boundary of a regular
neighborhood of the graph $\widehat G_s$ (this is checked inductively,
starting with the case $s=\emptyset$, and looking at the effect of
adding one edge at a time.)  Therefore their number equals the number
of connected components of the graph $\widehat G_s$, plus the rank of
the first homology of $\widehat G_s$, denoted $n(S)$. Then
$${\Phi}(G)\, =\, d^{V(G)/2}\, \sum_{s\subseteq E(\widehat G)} (-1)^{|s|} \, \frac{1}{d^{|s|}}\,\, d^{k(s)+n(s)}.$$

\begin{figure}[ht]
%\centering
\includegraphics[width=11cm]{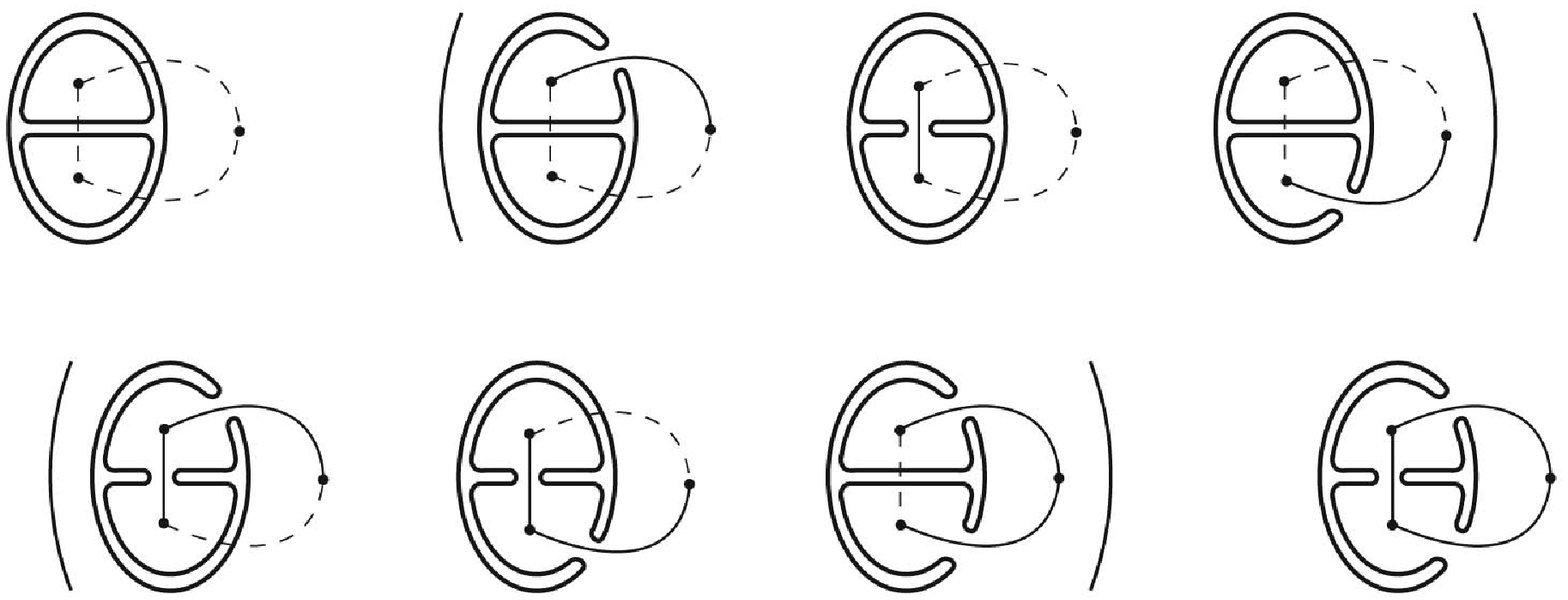}
\put(-325,91){$d$}
\put(-250,91){$-$}
\put(-161,91){$+$}
\put(-88,91){$+$}
\put(-330,22){$+\, \frac{1}{d}$}
\put(-240,22){$+$}
\put(-167,22){$+$}
\put(-80,22){$-\, \frac{1}{d^2}$}
     \caption{The expansions of $Q^{-1}\, {\chi}_Q(\widehat G)$, ${\Phi}(G)$ where $G$ is the theta graph in figure \ref{fig:ex5}.
     The edges of $\widehat G$ in each term which are in the given subset $s$ are drawn solid, other edges are dashed.}
\label{fig:ex6}
\end{figure}

We claim that the corresponding terms in the expansions of $Q^{-1}\, {\chi}_Q(\widehat G)$, ${\Phi}(G)$ are equal.
(See figure \ref{fig:ex6} for the expansions of $Q^{-1}\, {\chi}_Q(\widehat G)$, ${\Phi}(G)$ in the example $G=$ the theta graph
shown in figure \ref{fig:ex5}.) Since $G$ is a trivalent graph and we assumed $G$ is connected, its dual $\widehat G$ is
a triangulation (each face of $\widehat G$ has three edges), so that $2V(\widehat G)=F(\widehat G)+4=
V(G)+4$.
Therefore
$${\Phi}(G)\, =\, \sum_{s\subseteq E(\widehat G)} (-1)^{|s|} \, d^{V(\widehat G)-2+k(s)+n(s)-|s|}$$
The exponent simplifies because
%$$V(\widehat G)-2+k(s)+n(s)-|s| = 2k(s)-2.$$
$k(s)-n(s)+|s|=V(\widehat G)$.
This is obviously true for $s=\emptyset$, and inductively an addition of
one edge to $s$ either decreases $k(s)$ by $1$, or it increases $n(s)$
by $1$, so $k(s)-n(s)+|s|$ remains equal to $V(\widehat G)$. Because
$Q=d^2$,
$$Q^{-1}\, {\chi}_Q(\widehat G)\, =\, \sum_{s\subseteq E(\hat
G)} (-1)^{|s|} \, d^{2k(s)-2} = {\Phi}(G).$$
This concludes the proof of lemma \ref{commutative} for a connected graph $G$.

Now let $G$ be a trivalent, not necessarily connected, planar graph. To be specific, first suppose $G$
has two connected components, $G=G_1\sqcup G_2$. It follows from the definition of ${\Phi}$ that ${\Phi}(G)
={\Phi}(G_1)\cdot{\Phi}(G_2)$. Also note that $\widehat G$ is obtained from $\widehat G_1$, $\widehat G_2$
by identifying a single vertex. It is a basic property of the chromatic polynomial that in this situation
${\chi}_{\widehat G}(Q)=Q^{-1}{\chi}_{\widehat G_1}(Q)\cdot {\chi}_{\widehat G_2}(Q)$. Since the equality
(\ref{chromatic dual})
holds for connected graphs $G_1, G_2$, it also holds for $G$. This argument gives an inductive proof of
(\ref{chromatic dual}) for trivalent graphs with an arbitrary number of connected components.
\qed

\medskip

The {\em trace radical} of an algebra $A$ is the ideal consisting of
the elements $a$ of $A$ such that $tr(ab)=0$ for all $b\in A$.
% In the algebra $A$ (here $A={\mathcal C}_n^Q$ or $TL^d_n$), consider
The local relations on graphs which preserve the chromatic polynomial
of the dual, for a given value of $Q$, correspond to the elements of
the trace radical of $A={\mathcal C}_n^Q$ (see section \ref{Beraha}.)
It follows from lemma \ref{commutative} that the pullback by ${\Phi}$
of the trace radical in $TL^d$ is in the trace radical of ${\mathcal
C}^{d^2}$.

\begin{figure}[ht]
%\centering
\includegraphics[width=12cm]{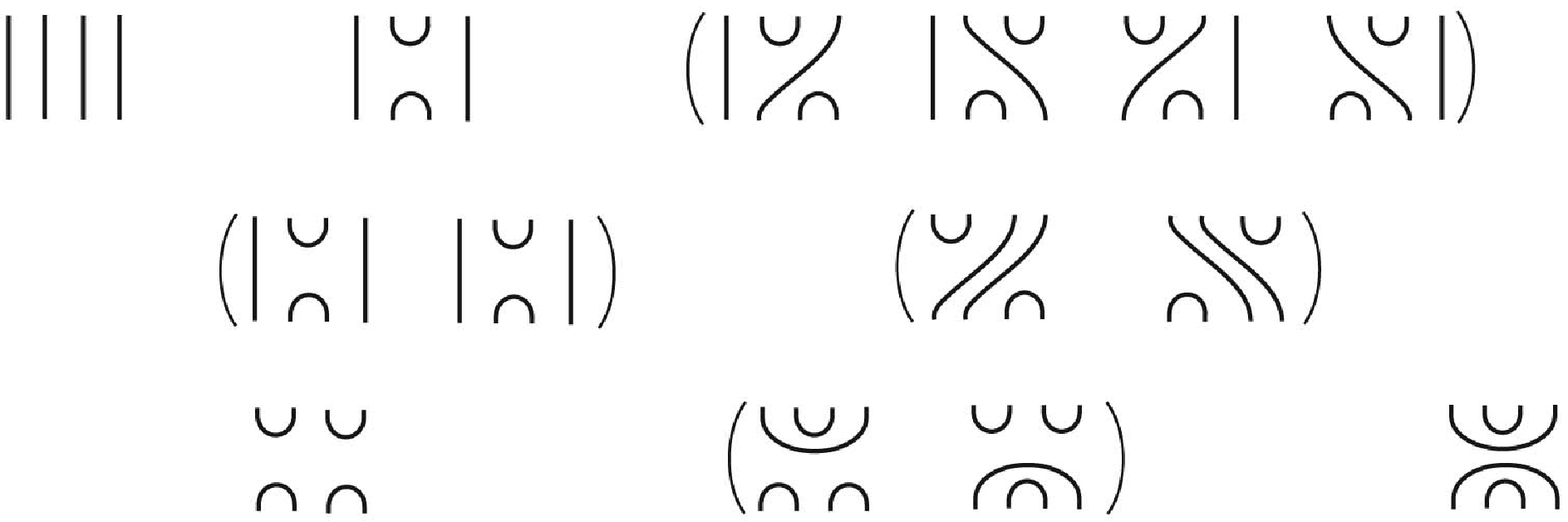}
{\large    \put(-383,96){$P^{(4)}\,=$}
    \put(-305,96){$-\,\frac{d}{d^2-2}$}
    \put(-230,96){$+\,\frac{1}{d^2-2}$}
    \put(-340,51){$+\,\frac{-d^2+1}{d^3-2d}$}
    \put(-198,51){$-\,\frac{1}{d^3-2d}$}
    \put(-350,10){$+\,\frac{d^2}{d^4-3d^2+2}$}
    \put(-250,10){$-\,\frac{d}{d^4-3d^2+2}$}
    \put(-90,10){$+\,\frac{1}{d^4-3d^2+2}$}}
     \caption{The Jones-Wenzl projector $P^{(4)}$ which generates the trace radical of $TL^{\phi}$. $P^{(4)}$ is defined
     for all values of $d$, and to get an element of $TL^{\phi}$ one sets $d={\phi}$ in the formula above. The
     homomorphism ${\Phi}\co {\mathcal C}^{{\phi}+1}\longrightarrow TL^{\phi}$ maps the relation
     (\ref{Tutte's identity4}): $\widehat Z_1\, +\, \widehat Z_2\; =\;
{\phi}^{-3}\,[\,\widehat Y_1\, +\, \widehat Y_2\,]$, dual to Tutte's identity (\ref{Tutte's identity3}), to the relation $P^{(4)}=0$.}
\label{fig:JW4}
\end{figure}

When $d={\phi}$, the trace radical of the Temperley-Lieb algebra is
generated by the Jones-Wenzl projector $P^{(4)}$ displayed in figure
\ref{fig:JW4}.  It is straightforward to check that $\Phi$ maps the
element (\ref{Tutte's identity4}) to $P^{(4)}$.  We have therefore
established Tutte's relation (\ref{Tutte's identity3}) as a
consequence of lemma \ref{commutative} and of the properties of the
Jones-Wenzl projector.  In section \ref{Beraha} we derive a recursive
formula for the identities of the chromatic polynomial at $Q$ obeying
(\ref{QJW}), which yields the identity (\ref{Tutte's identity3}) as a
special case.

\medskip

\section{Tutte's golden identity for the chromatic polynomial} \label{golden section}

In this section we give a proof of
Tutte's golden identity in the algebraic setting, established in the previous section. (Section
\ref{sec:levelrank} below presents an alternative proof using level-rank duality of the BMW algebras.)

\smallskip

\begin{theorem} \label{Tutte} \sl For a planar triangulation $\widehat G$,
\begin{equation}  {\chi}_{\widehat G}({\phi}+2)=({\phi}+2)\; {\phi}^{3\,V(\widehat G)-10}\; ({\chi}_{\widehat G}({\phi}+1))^2\end{equation}
where $V(\widehat G)$ is the number of vertices of $\widehat G$.
\end{theorem}

\smallskip

For the subsequent proof, it is convenient to reexpress this identity
in terms of the graph $G$ dual to $\widehat G$. Since
$\widehat G$ is a triangulation, $G$ is a connected trivalent graph.
Using the Euler characteristic, one observes that the number of
faces $F(\widehat
G)=2V(\widehat G)-4$. Since $V(G)=F(\widehat G)$,
the golden identity may be rewritten as
\begin{equation} \label{Tutte's identity2} {\chi}_{\widehat
    G}({\phi}+2)=
\frac{{\phi}+2}{\phi^4}\; {\phi}^{3 V(G)/2}\;
({\chi}_{\widehat G}({\phi}+1))^2
\end{equation}

\begin{figure}[ht]
%\centering
\includegraphics[height=1.7cm]{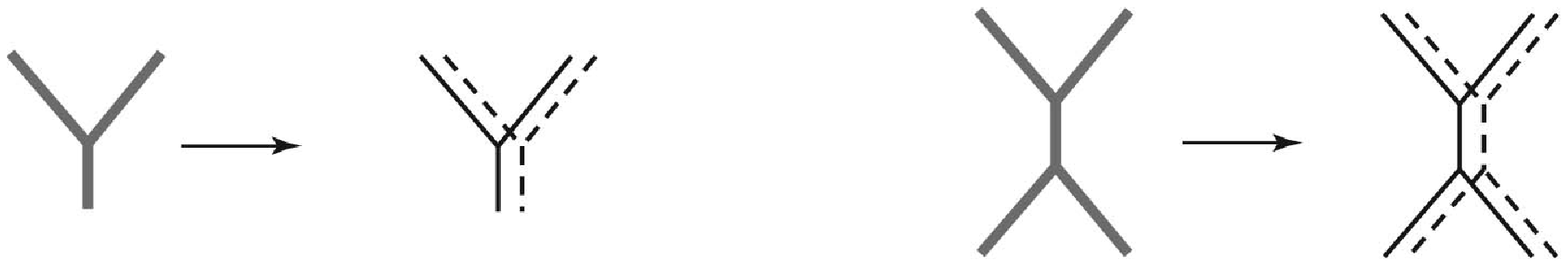}
{\scriptsize
    \put(-232,25){${\Phi}$}
    \put(-212,20){${\phi}^{3/2}\cdot$}
    \put(-62,26){${\Phi}$}
    \put(-40,20){${\phi}^3\cdot$}}
\caption{The map ${\Psi}\co {\mathcal FT}\longrightarrow {\mathcal FT}\times
{\mathcal FT}$. The different kinds of lines representing the graphs (gray, solid and dashed)
correspond to the three different copies of ${\mathcal FT}$.}
\label{fig:tutte}
\end{figure}

{\em Proof of theorem \ref{Tutte}.} Consider the vector space
${\mathcal FT}$ over ${\mathbb C}$ spanned by all connected planar
trivalent graphs. (There are no relations imposed among these graphs,
so this is an infinite dimensional vector space.)  Define a map $\Psi
\co {\mathcal FT}\longrightarrow {\mathcal FT}\times {\mathcal FT}$ on
the generators by \begin{equation} \label{Psi}{\Psi}(G)={\phi}^{3V/2}\; (G\times G),
\end{equation} where $V$
is the number of vertices of a trivalent graph $G$. (See figure \ref{fig:tutte}
illustrating the cases $V=1,2$.) Here ${\phi}$ denotes the golden
ratio.

\medskip

Consider the map ${\pi}_Q\co {\mathcal FT}\longrightarrow {\mathbb C}$
defined by taking the quotient of ${\mathcal FT}$ by the ideal
generated by the relations in the trivalent presentation of the
chromatic algebra, as given in definition \ref{chromatic definition}
and in figure \ref{fig:trivalent}. This quotient is $1-$dimensional
because the graphs in ${\mathcal FT}$ have no ends; applying the
relations allows any such graph to be reduced to a number.
Namely, this projection map ${\pi}_Q$ applied to a graph is the quantum
evaluation, or equivalently it is equal to $Q^{-1}$ times the
chromatic polynomial of the dual graph.

Tutte's identity follows from the following statement:

\begin{lem} \label{commutative1} \sl
The following diagram commutes
\[
\xymatrix{ {\mathcal FT}  \ar[d]^{\Psi} \ar[rr]^{\;\;\;{\pi}_{{\phi}+2}} && {\mathbb C}  \ar[d]^{=}\\
{\mathcal FT}\times {\mathcal FT}\ar[rr]^{\;\;\;({\pi}_{{\phi}+1})^2}  &&  {\mathbb C} }\]
\end{lem}

\medskip

{\em Proof:} The strategy is to check that the image under ${\Psi}$ of the three
relations in ${\mathcal FT}$ given in definition
\ref{chromatic definition} and figure \ref{fig:trivalent}
at $Q={{\phi}+2}$, hold
as a consequence of the relations in ${\mathcal FT}\times
{\mathcal FT}$ at $Q={\phi}+1$. Showing this means that $\frac{1}{{\phi}+2}\, \chi_{\widehat G}^{}(\phi+2)$ for any $G$ in ${\mathcal FT}$ can be evaluated instead in ${\mathcal FT} \times {\mathcal FT}$, and is equal to $$(({\phi}+1)^{-1}\; \chi_{\widehat G}^{}(\phi+1))^2=
\frac{1}{{\phi}^4}\, (\chi_{\widehat G}^{}(\phi+1))^2.$$
Checking two of the three relations is easy.
The value for $G=\,$circle at ${\phi}+2$ is ${\phi}+1$. The corresponding value
of ${\Psi}(G)=G\times G$ is ${\phi}^2={\phi}+1$. Likewise, tadpoles
vanish for any value of $Q$, so the image of the relation on the right
in figure \ref{fig:trivalent} under ${\Psi}$ clearly holds in
${\mathcal FT}\times {\mathcal FT}$.

\begin{figure}[ht]
\includegraphics[height=1.65cm]{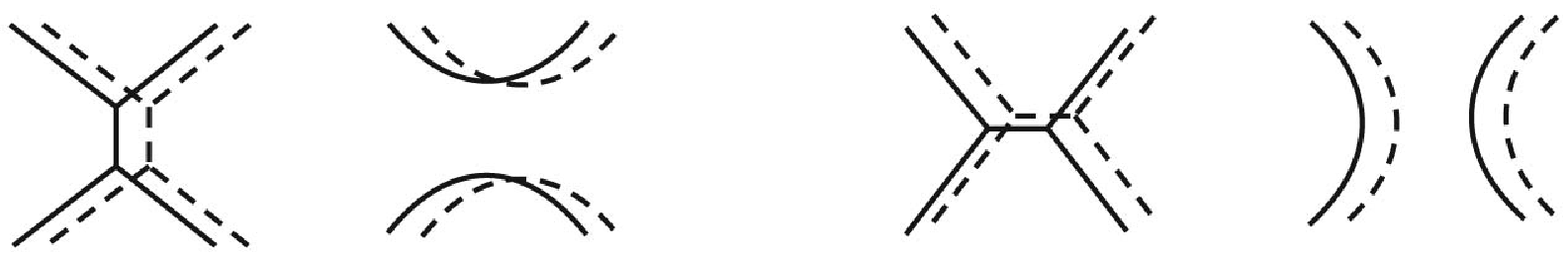}
    \put(-310,21){${\phi}^3\, \cdot$}
    \put(-241,21){$+$}
    \put(-166,21){$= \;{\phi}^3\, \cdot$}
    \put(-70,21){$+$}
\caption{The image under $\Psi$ of the H-I relation, defining the chromatic algebra (figure \ref{fig:trivalent}).}
\label{fig:tutte2}
\end{figure}

To derive the H-I relation, consider its image under $\Psi$, shown in figure \ref{fig:tutte2}.
Even though the image of ${\mathcal FT}$ under ${\pi}_Q$
is $1-$dimensional for any value of $Q$, we have seen that at
$Q={\phi}+1$ there is an {\em additional} relation (\ref{Tutte's
identity3}) obeyed by the chromatic polynomial.
(We will show in the next section that such additional relations exist
for any $Q$ satisfying (\ref{QJW}).)
This relation,
established in the previous section, is {\em not} a consequence of the
deletion-contraction rule, but is still consistent with the projection
map ${\pi}_{{\phi}+1}$. One checks using the first relation
in figure \ref{fig:trivalent} that (\ref{Tutte's identity4}) is
equivalent to each of the following two relations (in other words,
both of these relations are consistent with the projection
${\pi}_{{\phi}+1}$):
\begin{figure}[ht]
\includegraphics[height=1.6cm]{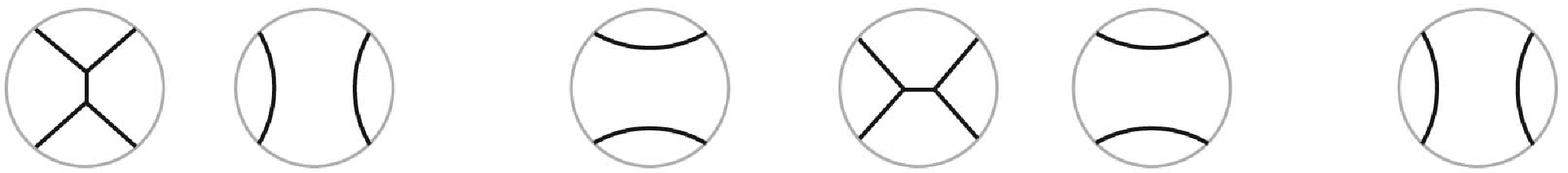}
{\small
    \put(-427,20){$\phi$}
    \put(-370,20){$=$}
    \put(-312,20){$+\,(1-{\phi})$}
    \put(-205,20){$\phi$}
    \put(-144,20){$=$}
    \put(-87,20){$+\,(1-{\phi})$}}
    \put(-223,7){$,$}
\caption{Relations in ${\mathcal C}_2^{{\phi}+1}$, equivalent to Tutte's relation (\ref{Tutte's identity3}).}
\label{fig:graphJW}
\end{figure}

Apply the relation on the left in figure \ref{fig:graphJW} to the expression in the left hand side in figure \ref{fig:tutte2},
in both
copies of ${\mathcal FT}$ at ${\phi}+1$, and use the identities involving the golden ratio $${\phi}+1={\phi}^2,\;\; {\phi}-1=
{\phi}^{-1}$$ to get

\begin{figure}[ht]
\includegraphics[height=1.55cm]{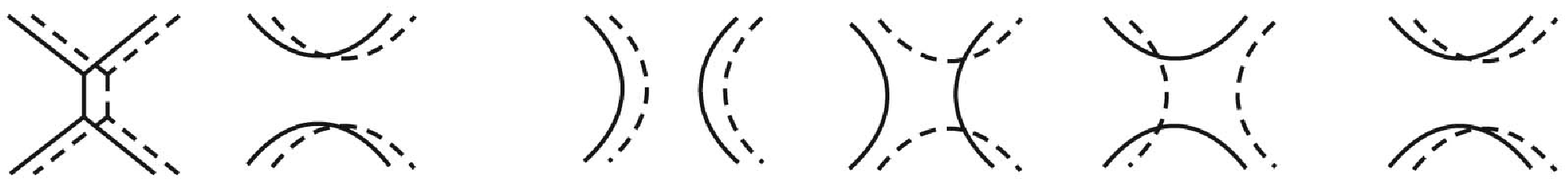}
    \put(-412,19){${\phi}^3\,\cdot$}
    \put(-348,19){$+$}
    \put(-283,19){$=\;{\phi}\,\cdot$}
    \put(-198,19){$-$}
    \put(-131,19){$-$}
    \put(-69,19){$+\;{\phi}\,\cdot$}
\caption{}
\label{fig:tutte3}
\end{figure}
The fact that the expression on the right is invariant under the 90 degree rotation establishes the
 relation in figure \ref{fig:tutte2}. This concludes the proof of lemma \ref{commutative1} and of theorem \ref{Tutte}.
\qed

{\bf Remarks}. 1. The proof of the theorem above could also have been given directly
in the context of algebras and their traces. The proof of lemma
\ref{commutative1} above shows the existence of a map ${\mathcal
  C}^{{\phi}+2}\longrightarrow ({\mathcal C}^{{\phi}+1}/R)\times ({\mathcal
  C}^{{\phi}+1}/R)$, where $R$ denotes the trace radical. One may consider then a diagram involving the
Temperley-Lieb algebras, analogous to (\ref{chromaticTL}), and theorem
\ref{Tutte} follows from applying the trace to these algebras.

2. The map ${\Psi}$, considered in the proof of the theorem \ref{Tutte}, has a simple
definition (\ref{Psi}) only in the context of trivalent graphs. The chromatic algebra
may be defined in terms of all, not just trivalent, graphs (see the next section), however the
extension of ${\Psi}$ to vertices of higher valence is substantially more involved.

\bigskip

\section{Level-rank duality}
\label{sec:levelrank}

In the previous section we gave an algebraic proof of Tutte's identity
(\ref{Tutte's identity1}) using the chromatic algebra, without making
reference to the Birman-Murakami-Wenzl algebra. In this section we
explain how an essential ingredient for a conceptual understanding of
this identity is the {\em level-rank duality} of the $SO(N)$
topological quantum field theories and of the BMW algebras.
The golden identity (\ref{Tutte's identity1}) is a
consequence of the duality between the $SO(3)_4$ and $SO(4)_3$ algebras.
This viewpoint explains why such an identity relates the
chromatic polynomial ${\chi}(Q)$ at the values $Q={\phi}+1, {\phi}+2$,
and why one does not expect a generalization for other values of
$Q$.

We will work with $SO(N)$
Birman-Murakami-Wenzl algebras which underlie the construction of the
(doubled) $SO(N)$ TQFTs (see for example \cite{F} for a related
discussion concerning the TL algebra and the $SU(2)$ theories, and
\cite[section 3.1]{FK} in the case of $BMW(3)$.) Note that the
chromatic algebra and the $SO(3)$ BMW algebra are closely related
\cite{FK}.

\begin{figure}[ht]
%\centering
\includegraphics[height=1.05cm]{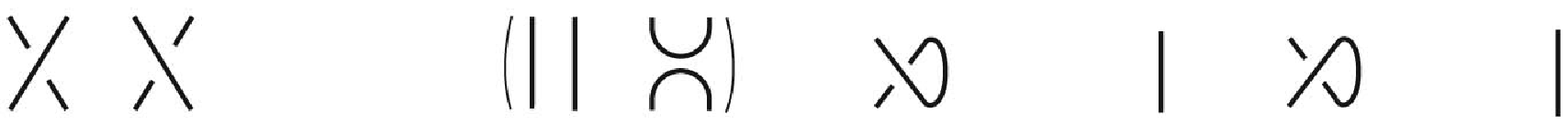}
    \put(-377,12){$-$}
    \put(-338,12){$=\; (q-q^{-1})$}
    \put(-246,12){$-$}
    \put(-207,5){,}
    \put(-148,12){$= q^{1-N}$}
    \put(-99,5){,}
    \put(-46,12){$= q^{N-1}$}
\caption{}
\label{fig:BMWskein}
\end{figure}

The $SO(N)$ BMW algebra
$BMW(N)_n$ is the algebra of framed tangles on $n$
strands in $D^2\times [0,1]$ modulo regular isotopy and the $SO(N)$ Kauffman
skein relations in figure \ref{fig:BMWskein} \cite{BW,M}.
By a tangle we mean a collection of curves (some of them perhaps
closed) embedded in $D^2\times [0,1]$, with precisely $2n$ endpoints,
$n$ in $D^2\times\{0\}$ and $D^2\times\{1\}$ each, at the prescribed
marked points in the disk. The tangles are framed, i.e.\ they are given with a
trivialization of their normal bundle. (This is necessary since the last
two relations in figure \ref{fig:BMWskein} are not invariant under the first
Reidemeister move.) As with $TL$, the multiplication is given
by vertical stacking, and we set $BMW(N)=\cup_n BMW(N)_n$.

The trace, $\hbox{tr}_K\co BMW(N)_n\longrightarrow {\mathbb C}$,
is defined on the generators
(framed tangles) by connecting the top and bottom endpoints by
standard arcs in the complement of $D^2\times[0,1]$ in $3$-space,
sweeping from top to bottom, and computing the $SO(N)$ Kauffman
polynomial (given by the skein relations above) of the resulting
link. The skein relations imply that
deleting a circle has
the effect of multiplying the element of $BMW(N)$ by
\begin{equation}
d_{N}=1+\frac{q^{N-1}-q^{-(N-1)}}{q-q^{-1}}\ .
\label{dN}
\end{equation}
In fact, it follows from the first relation in figure
\ref{fig:BMWskein} that
the subalgebra
of $BMW(N)$ generated by
$e_i = 1- (B_i-B^{-1}_i)/(q-q^{-1})$ is
TL algebra $TL^{d_N}$.

A planar version of the $SO(N)$ BMW algebra, generated by
$4$-valent planar graphs, is found by using
\begin{equation} \label{basischange}(\crossing)\; =\;
 q(\smoothing)-(\slashoverback)
+q^{-1}(\hsmoothing)\; = \;
 q^{-1}(\smoothing)-(\backoverslash)
+q(\hsmoothing)
\end{equation}
When $N=3$, the relations in this
planar version reduce to the chromatic algebra relations \cite{FR,FK}.

%; for example, the
%trace radical is non-trivial, and for $N=3$ the Jones-Wenzl projectors
%pulled back from the Temperley-Lieb algebra are in the trace radical.

When $q$ is a root of unity, the BMW algebra has a variety of
special properties. A key property for us occurs when the ``level''
$k$, defined via
$$q=e^{i\pi/(k+N-2)},$$
is an integer. The algebra
$BMW(N)$ for $k$ integer (which we label by $BMW(N,k)$)
is related to $SO(N)$ Chern-Simons topological field
theory at level $k$ and the
Wess-Zumino-Witten conformal field theory with symmetry algebra given
by the Kac-Moody algebra $SO(N)_k$ \cite{Witten,W}.
Note that with this definition of $k$, we have $q^{N-1}=- q^{-(k-1)}$.
It is then easy to check if a braid $B$ obeys the relations in
$BMW(N,k)$, then $-B^{-1}$ obeys $BMW(k,N)$.

We use this special property to define an algebra homomorphism
${\psi}\co BMW(N,k)\longrightarrow BMW(k,N)$, which sends an additive
generator (a framed tangle in $D^2\times [0,1]$) to a tangle with each
crossing reversed, and multiplied by $(-1)^{\# {\rm crossings}}$.
\begin{lem} \label{lem:level rank}\sl
The map $\psi$ defined above is an isomorphism, and moreover it preserves the trace:
\begin{equation} \label{level rank diagram} \xymatrix{ BMW(N,k)  \ar[d]^{tr} \ar[r]^{\psi} & BMW(k,N)  \ar[d]^{tr}\\
{\mathbb C}\ar[r]^=  &  {\mathbb C} }
\end{equation}
\end{lem}

Recall that the traces of the two algebras are the $SO(N)$,
respectively $SO(k)$, Kauffman polynomials evaluated at
$q=e^{i\pi/(k+N-2)}$. The proof of the lemma is checked directly from
definitions: ${\psi}$ is well-defined with respect to the relations in
figure \ref{fig:BMWskein} and the last two Reidemeister moves, and
${\psi}$ has the obvious inverse. The diagram above commutes since
computing the trace (the Kauffman polynomial) involves the relations
in figure \ref{fig:BMWskein}.  Note that the weight of a circle $d_N$
from (\ref{dN}) written in terms of $k$ does not change if $N$ and $k$
are interchanged. This isomorphism is an algebra analogue of {\em
level-rank duality} (for $SO(N)$ TQFTs).

The golden identity (\ref{Tutte's identity1}) comes from studying the
very special case of the duality between $SO(3)_4$ and $SO(4)_3$
theories, where $q=\exp(i\pi/5)$, and then exploiting a known
but possibly underappreciated fact about the $SO(4)$
BMW algebra: it can be decomposed into the product of two TL
algebras. This decomposition generalizes
the isomorphism between the Lie algebras
$SO(4)$ and $SU(2)\times SU(2)$, where
the fundamental four-dimensional representation of $SO(4)$ decomposes
into the product of two spin-1/2 representations of $SU(2)$.
The homomorphism $\Psi$ used in the previous section is the chromatic
analogue of the composite isomorphism from $SO(3)_4$ to
$SO(3)_{3/2}\times SO(3)_{3/2}$.

\begin{figure}[ht]
%\centering
\includegraphics[height=1cm]{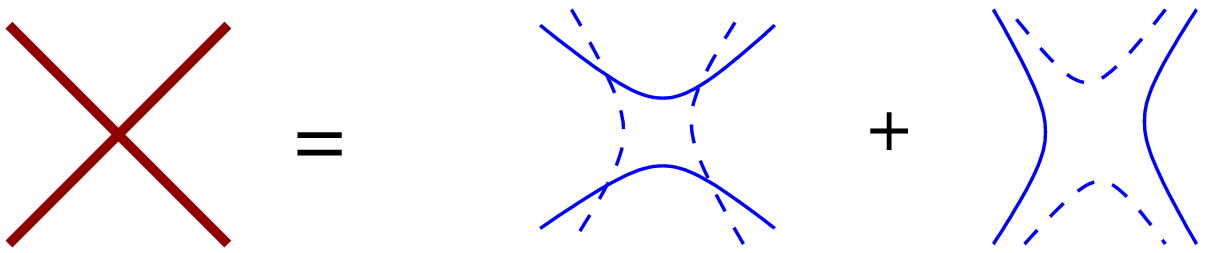}
\caption{The crossing in $BMW(4)$ mapped to $TL\times TL$}
\label{fig:doubled}
\end{figure}

To make this precise, consider the
``skein'' presentation, $TL^{\rm skein}$, of the Temperley-Lieb
algebra, given by framed tangles in $D^2\times [0,1]$ modulo regular
isotopy and the skein relations defining the Kauffman bracket (see
section 2 in \cite{FK}).  The map $BMW(4)^q_n \to TL^d_n \times
TL^d_n$ takes each generator $g$, $g=B,B^{-1}$ or $e$ to the product
of the corresponding elements in $TL^{\rm skein}$,
$g\mapsto(g,g)$. The parameters $q$ and $d$ in these algebras are
related by $d=-q-q^{-1}$ as before; note that the weight
$d_4=(q+q^{-1})^2$ of a circle in $BMW(4)$ from (\ref{dN}) is indeed
the square of the weight of a circle in $TL^d$. This
map is well-defined: the Reidemeister moves and the BMW
relations in figure \ref{fig:BMWskein} hold in $TL\times TL$ as a
consequence of the Reidemeister moves and and the Kauffman bracket
skein relations.  The same argument shows that this map preserves the
trace.  In the planar formulation, the crossing defined by
(\ref{basischange}) is expressed as a linear combination in $TL\times
TL$, as displayed in figure \ref{fig:doubled}.
We denote one
of the TL copies with a dashed line to emphasize that these copies are
independent.

%The combination of level-rank duality with this map of $BMW(4)$ to
%$TL\times TL$ implies that the algebra $BMW(N,4)$ for any $N$ can be
%expressed in terms of doubled strands. In the planar version, the
%crossing in $BMW(N,4)$ is then rewritten in terms of the doubled
%lines, akin to the relation for $BMW(4)$ in figure
%\ref{fig:doubled}. (The $BMW(N,4)$ case is more involved since one must map
%$B\to -B^{-1}$ in the level-rank duality.) However, this is quite
%useful, because solutions of the Yang-Baxter equation follow from
%representations of the BMW algebra \cite{Wadati}. Those corresponding
%to $BMW(4,k)$ correspond to coupled Potts or coupled loop models, and
%using level-rank duality allows one to find an integrable critical
%theory \cite{FJ}.

\begin{figure}[ht]
\includegraphics[height=1.6cm]{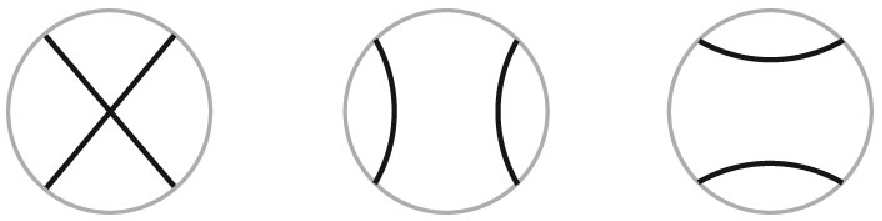}
{\small
    \put(-193,20){$\phi$}
    \put(-130,20){$=$}
    \put(-60,20){$+$}}
\caption{}
\label{fig:graphJW1}
\end{figure}

The final ingredient in this proof of Tutte's identity (\ref{Tutte's
identity1}) is the isomorphism $BMW(3,3/2)/R\cong TL$, where $R$
denotes the trace radical. Using the relation between the $SO(3)$ BMW
algebra and the Temperley-Lieb algebra \cite{FK}, one observes that
the relation in figure \ref{fig:graphJW1} is in the trace radical of
the (planar presentation of) $BMW(3,3/2)$. (In fact, this relation is
equivalent to Tutte's relation (\ref{Tutte's identity3}).)  Using this
relation, $BMW(3,3/2)/R$ is mapped to $TL$ by resolving all $4$-valent
vertices. The inverse map, $TL\longrightarrow BMW(3,3/2)$, is given by
mapping the generators $e_i$ of $TL$ to the same generators of the BMW
algebra.

The golden identity follows from the duality between $BMW(3,4)$ and
$BMW(4,3)$, and then mapping the latter into a product of two copies
of $BMW(3,3/2)/R$. The map from $BMW(4)^q$ to $TL^d\times TL^d$ above
exists for any value of $q$, and at the special value $q=\exp(i\pi/5)$
we showed above that this map lifts to a homomorphism
$$BMW(4,3)\longrightarrow (BMW(3,3/2)/R)\times (BMW(3,3/2)/R).$$ The
trace of the $SO(3)$ BMW algebra is given (up to a normalization) by
the chromatic polynomial, see theorem 6.3 and corollary 6.5 in
\cite{FK}.  Tutte's identity then follows by applying the trace to the
algebras above: the trace of $BMW(4,3)\cong BMW(3,4)$ is given by the
chromatic polynomial evaluated at $Q=\phi+2$, while the product of the
traces on the right corresponds to $({\chi}({\phi}+1))^2$.

The fact that the golden identity relates a chromatic polynomial
squared to the chromatic polynomial at another value is very special
to these values of $Q$.  Generalizations using level-rank duality
$SO(N)_4\to SO(4)_N$ for other values of $N$ do not involve the
chromatic polynomial.  Generalizations using $SO(N)_3\to SO(3)_N$ at
other $N$ give a linear relation between the chromatic polynomial and
the Markov trace for the $SO(N)_3$ BMW algebra, very unlike the golden
identity.

\bigskip

\section{Tutte's relations and Beraha numbers} \label{Beraha}

In this section we establish relations for the chromatic polynomial
evaluated at any value of
\begin{equation}  \label{generalizedBeraha}
Q=2+2\cos\left(\frac{2\pi j}{n+1}\right),
\end{equation}
for any positive integers obeying $j< n$. Each such relation is
independent of (but consistent with) the contraction-deletion
relation, and Lemmas \ref{recursive}, \ref{recursive2} give a
recursive formula for them. Our result generalizes Tutte's identity
(\ref{Tutte's identity3}) to this set of $Q$; note that the values
(\ref{generalizedBeraha}) are dense in the interval $[0,4]$ in the
real line. When $j=1$, these are the Beraha numbers, which strong
numerical evidence suggests are the accumulation points of the zeros
of $\chi(Q)$ for planar triangulations.  Tutte had conjectured such a
relation would exist for all Beraha numbers, and found examples in
several cases \cite{T3}. Our results follow from the observation that
the chromatic polynomial relations for a given value of $Q$ correspond
to elements of the trace radical of the chromatic algebra ${\mathcal
  C}^Q$.  To derive these relations, it is convenient to define the
chromatic algebra in terms of all planar graphs, not just the
trivalent ones used above.

\subsection{A presentation of the chromatic algebra using the
  contraction-deletion rule.}
Consider the set ${\mathcal G}_n$ of the isotopy classes
of planar graphs $G$ embedded in the rectangle $R$ with $n$ endpoints
at the top and $n$ endpoints at the bottom of the rectangle. (The
intersection of $G$ with the boundary of $R$ consists precisely of
these $2n$ points, figure 3, and the isotopy of graphs is required to preserve the boundary.)
It is convenient to divide the set of
edges of $G$ into {\em outer} edges, i.e. those edges that have an
endpoint on the boundary of $R$, and {\em inner} edges, whose vertices
are in the interior of $R$.

The relations in \ref{chromatic definition} defining the chromatic
algebra apply to trivalent graphs. Instead of generalizing them
directly, we define the chromatic algebra here using the
contraction-deletion rule. It is shown in \cite{FK} that the two
definitions (\ref{chromatic definition}, \ref{chromatic definition1})
give isomorphic algebras. Analogously to the other definition, the
idea is to view the contraction-deletion rule (\ref{chromatic poly1})
as a linear relation between the graphs $G, G/e$ and $G\backslash e$.
To make this precise, let ${\mathcal F}_n$ denote the free algebra
over ${\mathbb C}[Q]$ with free additive generators given by the
elements of ${\mathcal G}_n$.  The multiplication is given by vertical
stacking. Define ${\mathcal F}=\cup_n {\mathcal F}_n$.  Consider the
following set of local relations on the elements of ${\mathcal G}_n$.
(Note that these relations only apply to {\em inner} edges which do
not connect to the top and the bottom of the rectangle.)

\begin{figure}[ht]
\includegraphics[height=2cm]{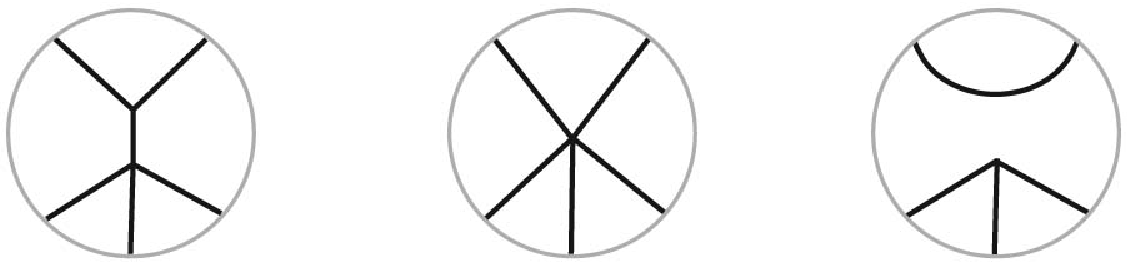}
    \put(-175,28){$=$}
    \put(-80,28){$-$}
{\scriptsize
    \put(-212,26){$e$}
    \put(-245,54){$G$}
    \put(-152,54){$G/e$}
    \put(-7,54){$G\backslash e$}}
\caption{Relation (1) in the chromatic algebra}
\label{fig:chromedef1}
\end{figure}

\medskip

(1) If $e$ is an inner edge of a graph $G$ which is not a loop, then
$G=G/e-G\backslash e$, as illustrated in figure \ref{fig:chromedef1}.

\begin{figure}[ht]
\includegraphics[height=2cm]{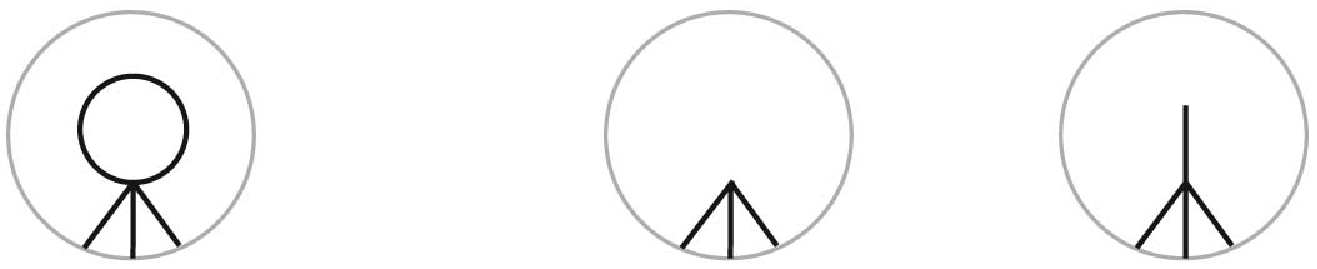}
    \put(3,25){$=\; 0.$}
    \put(-218,25){$=\; (Q-1)\;\cdot$}
    \put(-100,15){$,$}
\caption{Relations (2), (3) in the chromatic algebra}
\label{fig:chromedef23}
\end{figure}

(2) If $G$ contains an inner edge $e$ which is a loop, then $G=(Q-1)\;
G\backslash e$, as in figure \ref{fig:chromedef23}. (In particular,
this relation applies if $e$ is a loop not connected to
the rest of the graph.)

(3) If $G$ contains a $1$-valent vertex (in the interior of the rectangle)
as in figure \ref{fig:chromedef23}, then $G=0$.

\medskip

\begin{defi} \label{chromatic definition1}
The {\em chromatic algebra} in degree $n$, $\overline{\mathcal C}_n$, is an algebra over ${\mathbb C}[Q]$
which is defined as the quotient of the free algebra ${\mathcal F}_n$ by the ideal
$I_n$ generated by the relations (1), (2), (3). $\bar{\mathcal C}^Q_n$ denotes the algebra over ${\mathbb C}$
when $Q$ is specialized to a complex number. Set $\overline{\mathcal C}=\cup_n \overline{\mathcal C}_n$.
Analogously to section \ref{sec:TL}, the trace, $tr_{\chi}\co \bar{\mathcal C}^Q\longrightarrow{\mathbb C}$ is defined on the additive
generators (graphs) $G$ by connecting the
endpoints of $G$ by arcs in the plane (denote the result by $\overline G$)
and evaluating $$Q^{-1}\cdot {\chi}_{\widehat{\overline G}}(Q).$$
\end{defi}

One checks that the trace is well-defined with respect to the relations $(1)-(3)$. For example, the relation (1) corresponds
to the contraction-deletion rule for the chromatic polynomial of the dual graph: ${\chi}_{\widehat G}=
{\chi}_{\widehat G\backslash \widehat e} - {\chi}_{\widehat G/\widehat e}$, where $\widehat e$ is the edge of $\widehat G$
dual to $e$. The relation $(2)$ corresponds to deleting a $1$-valent vertex and the adjacent edge of the dual
graph, and the chromatic polynomial vanishes in case $(3)$ since the dual graph has a loop.

{\bf Remark.} The trace may also be described in terms of the {\em flow polynomial} of $\overline G$. Both the flow polynomial
and the chromatic polynomial are one-variable specializations of the two-variable Tutte polynomial \cite{B}.

\begin{figure}[ht]
\includegraphics[width=12cm]{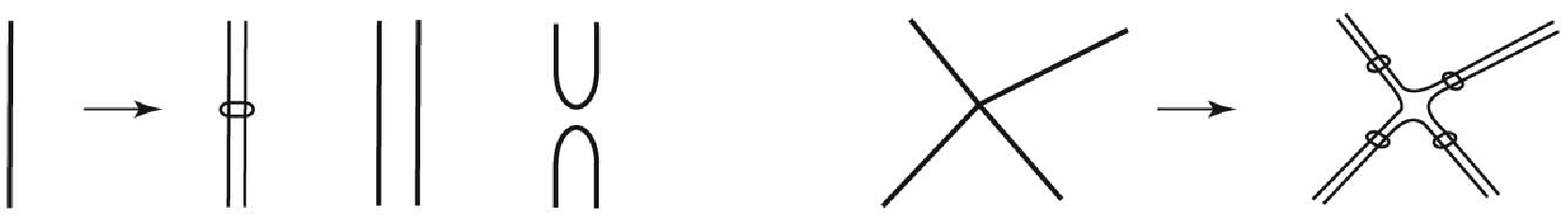}
    \put(-322,29){$\overline{\Phi}$}
    \put(-277,21){$=$}
    \put(-244,21){$-\frac{1}{d}$}
    \put(-85,29){$\overline{\Phi}$}
    \put(-62,21){$d\;\cdot$}
\caption{Definition of the homomorphism $\overline{\Phi}\co \bar{\mathcal C}^Q_n\longrightarrow TL^{\sqrt Q}_{2n}$}
\label{fig:anotherphi}
\end{figure}

Consider the algebra homomorphism $\overline{\Phi}\co \bar{\mathcal
  C}^{d^2}_n\longrightarrow TL^d_{2n}$, analogous to the one given in
section \ref{sec:TL} for trivalent graphs shown in figure \ref{fig:anotherphi}.
The factor in the definition of $\overline{\Phi}$ corresponding to a
$k$-valent vertex is $d^{(k-2)/2}$, for example it equals $d$ for the
$4$-valent vertex in figure \ref{fig:anotherphi}. The overall factor
for a graph $G$ is the product of the factors $d^{(k(V)-2)/2}$ over
all vertices $V$ of $G$.

One checks that $\overline{\Phi}$ is well-defined. For example, for the
defining relation (1) of the chromatic algebra in figure
\ref{fig:chromedef1}, one applies $\overline{\Phi}$ to both sides and expands
the projector at the edge $e$, as shown in figure
\ref{fig:anotherphi}. The resulting relation holds due to the choices
of the powers of $d$ corresponding to the valencies of the
vertices. Similarly, one checks the relations (2) and (3).

It is shown in \cite[Section 4]{FK} that the inclusion $\{$trivalent
graphs$\}\subset \{$all graphs$\}$ induces an isomorphism of the
algebras ${\mathcal C}^Q_n$, $\bar {\mathcal C}^Q_n$ in definitions \ref{chromatic definition},
\ref{chromatic
definition1}. Therefore from now on we will often use the same notation, ${\mathcal C}^Q_n$, for
both. It follows that the version of lemma \ref{commutative}
holds for the algebra $\bar{\mathcal C}^Q_n$ defined in \ref{chromatic
definition1}. Specifically, in the notations of lemma
\ref{commutative}, for {\em any} planar graph $G$, $Q^{-1}\,
{\chi}_Q(\widehat G)\, =\, \overline{\Phi}(G)$, and therefore $\overline{\Phi}$
preserves the algebra traces:
\begin{equation} \label{another commutative}
\xymatrix{ \bar{\mathcal C}_n^{\scriptsize Q}  \ar[d]^{tr_{\chi}} \ar[r]^{\overline{\Phi}} & TL_{2n}^{\sqrt Q}  \ar[d]^{tr_{d}}\\
{\mathbb C}\ar[r]^=  &  {\mathbb C} }\end{equation}
This is true since any graph generator of $\bar {\mathcal C}^Q_n$ in definition \ref{chromatic definition1} is equivalent
to a linear combination of {\em trivalent} graphs, using the contraction-deletion rule (1). The commutativity of the
diagram above then follows from lemma \ref{commutative} and the fact that $\overline{\Phi}$ and the traces ${tr_{\chi}}$,
$tr_{d}$ of the two algebras are well-defined.

\subsection{Chromatic polynomial relations and the trace radical.}
It follows from (\ref{another commutative}) that the pull-back of the trace radical in $TL_{2n}^{\sqrt Q}$ is in the
trace radical of ${\mathcal C}_n^{\scriptsize Q}$. (Recall that the trace radical in the algebra $A$,
where $A={\mathcal C}_n^Q$ or $TL^d_n$,
is the ideal consisting of the elements $a$ of $A$ such that $tr(ab)=0$ for all $b\in A$.)

Next observe that the local relations on graphs which preserve the
chromatic polynomial of the dual, for a given value of $Q$, correspond
to the elements of the trace radical in ${\mathcal C}^Q$. (By
dualizing the relation, one gets relations which preserve the
chromatic polynomial of the graphs themselves, as opposed to that of
their dual graphs; for example see figures 1 and 2.)  Indeed, suppose
$R=\sum_i a_i G_i$ is a relation among graphs in a disk $D$, so each
graph $G_i$ has the same number of edges meeting the boundary of the
disk. Suppose first that this number is even, say equal to $2n$.
Divide the boundary circle of the disk into two intervals, so that
each of them contains precisely $n$ endpoints of the edges.  Consider
the disk $D$ as a subset of the $2-$sphere $S^2$.  The fact that $R$
is a relation means that for any fixed graph $G$ in the complement
$S^2\smallsetminus D$, with the same $2n$ points on the
boundary, the linear combination $\sum_i a_i {\chi}_{\widehat{G_i\cup
G}}$ vanishes at $Q$.  Since both the disk $D$ and its complement
$S^2\smallsetminus D$ are homeomorphic to a rectangle, one may
consider $R$ as an element in ${\mathcal C}^Q_n$, and moreover it is
in the trace radical: $tr(R\cdot G)=0$, and the graphs $G$ additively
generate ${\mathcal C}_n$.

Converting the disk $D$ above into a rectangle, the subdivision of the
$2n$ boundary points into two subsets of $n$, and the fact that this
number is even, may seem somewhat artificial. This reflects the
algebraic structure of the setting we are working in.  The discussion
may be carried through in the context of the chromatic {\em category},
and further {\em planar algebra}, where the algebraic structure is
more flexible while the notion of the trace radical is retained. For
example, in the category the multiplication (vertical stacking) is
complemented by tensor multiplication (horizontal stacking).
Describing these structures in further detail would take us outside
the scope of the present paper, so instead we refer the interested
reader to \cite[section 2]{F} where the discussion is given in the
similar context of the Temperley-Lieb algebra.  One observes that a
relation preserving the chromatic polynomial, in fact, corresponds to
an element of the ideal closure of the trace radical in the chromatic
category, not just algebra. (The importance of this distinction will
become clear in the following subsections.)  The converse argument
shows that an element of the trace radical may be viewed as a relation
among planar graphs, preserving the chromatic polynomial of the duals.

\subsection{The trace radical in $TL_{2m}^d$ and relations in ${\mathcal
    C}^{Q}_m$.}
The structure of the trace radical in the Temperley-Lieb algebra is
well-understood. In particular, they occur for each special value of
$d$ defined by
\begin{equation} \label{special d}
d\, = \, 2\, \cos \left(\frac{\pi j}{n+1}\right),
%\, =\, ({2+2\,\cos(2j{\pi}/(n+1))})^{1/2},
\end{equation}
where $j$ and $n$ are positive integers obeying $j<n$.
When $j$ and $n+1$ are coprime, the trace radical in $TL^d$
is generated by an element $P^{(n)}$ called the {\em Jones-Wenzl
projector} \cite{Jo,We}. A theorem of Goodman-Wenzl
\cite[Appendix]{F} shows that for values of $d$ other than (\ref{special
d}), the Temperley-Lieb category does not have any non-trivial proper
ideals. We briefly review the basic properties of these projectors
below; see \cite{KL} for more details, and \cite{F} for a discussion
of the trace radical in the context of the Temperley-Lieb category.

\begin{figure}[ht]
\includegraphics[width=6cm]{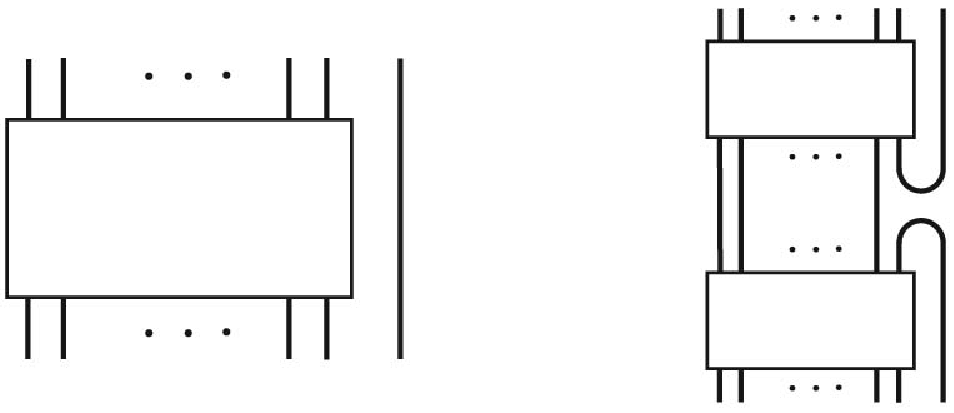}
    \put(-212,31){$P^{(n)}\, =$}
    \put(-152,31){$P^{(n-1)}$}
    \put(-87,31){$-\; \frac{{\Delta}_{n-2}}{{\Delta}_{n-1}}$}
{\small
    \put(-40,54){$P^{(n-1)}$}
    \put(-40,12){$P^{(n-1)}$}}
     \caption{A recursive formula for the Jones-Wenzl projectors in the Temperley-Lieb algebra.}
\label{fig:JonesWenzl}
\end{figure}

The Jones-Wenzl projector $P^{(r)}_i$ acts in $TL_n$ for any value of
$d$ with $n\ge r$; when $n>r$ we include the subscript to indicate
that it is acting non-trivially on the strands labeled $i,i+1,\dots
i+r-1$. The first two Jones-Wenzl projectors are
$P^{(1)}=1$, and $P_i^{(2)}=1-E_i$. (We use the notation for the
generators of the Temperley-Lieb algebra introduced in section \ref{TL subsection}.)
A recursive
formula (cf \cite{KL}) giving the rest is
\begin{equation}
P^{(n)} = P^{(n-1)}_1 -\,  \frac{d\,  \Delta_{n-2}}{\Delta_{n-1}}\,
P^{(n-1)}_1\,  E_{n-1}\,  P^{(n-1)}_1,
\label{JonesWenzl}
\end{equation}
where the number $\Delta_n$ is simply the trace:
\begin{equation}
\Delta_n\, =\, \hbox{tr}_d \; P^{(n)}\ .
\end{equation}
This recursion relation is illustrated
in figure \ref{fig:JonesWenzl}.
Taking the trace of the elements in the recursion relation yields
$\Delta_n=d\Delta_{n-1} - \Delta_{n-2}$.
Since $\Delta_1=d$ and $\Delta_2=d^2-1$,
\begin{equation}
\Delta_n=\frac{\sin[(n+1)\theta]}{\sin(\theta)}
\label{Deltatheta}
\end{equation}
where $\theta$ is defined via $d=2\cos(\theta)$.

Of course any element in the trace radical has its own trace equal to zero.
From
the explicit formula (\ref{Deltatheta}) for $\Delta_n$, it follows
that $P^{(n)}$ can be in the trace radical only when $\theta=\pi
j/(n+1)$, where $j$ is a non-zero integer not a multiple of $n+1$. In
terms of $d$, these correspond to the values in $(\ref{special d})$.
Indeed, for these values of $d$, $P^{(n)}$ generates the
trace radical, and moreover the theorem of Goodman-Wenzl guarantees that
this is the
unique proper ideal in the category for coprime $j$ and $n+1$.

{\bf Remark.} A useful way of understanding the Jones-Wenzl projectors
is to think of each strand in the Temperley-Lieb algebra as carrying
spin-1/2 of the quantum-group algebra $U_q(sl(2))$, where
$d=q+q^{-1}$.  There exist spin-$r/2$ representations of this algebra
behaving similarly to those of ordinary $su(2)$, except for the fact
at the special values of $d$ in (\ref{special d}) they are irreducible
only for $r<n$. In the algebraic language,
the recursion relation describes taking a tensor product of
representations, so that $P^{(r)}$ is the projector onto the
largest-possible value of spin $r/2$ possible for $r$ strands. For
example, $P^{(2)}_i$ projects the strands $i,i+1$ onto spin $1$, while
the orthogonal projector $E_i$ projects these two strands onto spin
$0$. The fact that the spin-$n/2$ representation at a
special value of $d$ is reducible is the reason the corresponding
projector $P^{(n)}$ can be set to zero.
The strands in the chromatic algebra can be viewed as
carrying spin 1. It is logical to expect that a projector from $m$
spin-1 strands onto spin $m$ exists, and that it generates a proper
ideal. We will show that this is indeed so in the following.

\begin{lem} \label{imagelemma}
\sl The Jones-Wenzl projector $P^{(2m)}\in TL_{2m}$ is in the image of the chromatic algebra: $P^{(2m)}\in {\Phi}({\mathcal C}_m)$,
so its pullback at the corresponding special value of $Q$ is in the trace
radical of the chromatic algebra.
\end{lem}

\begin{figure}[ht]
\includegraphics[height=2cm]{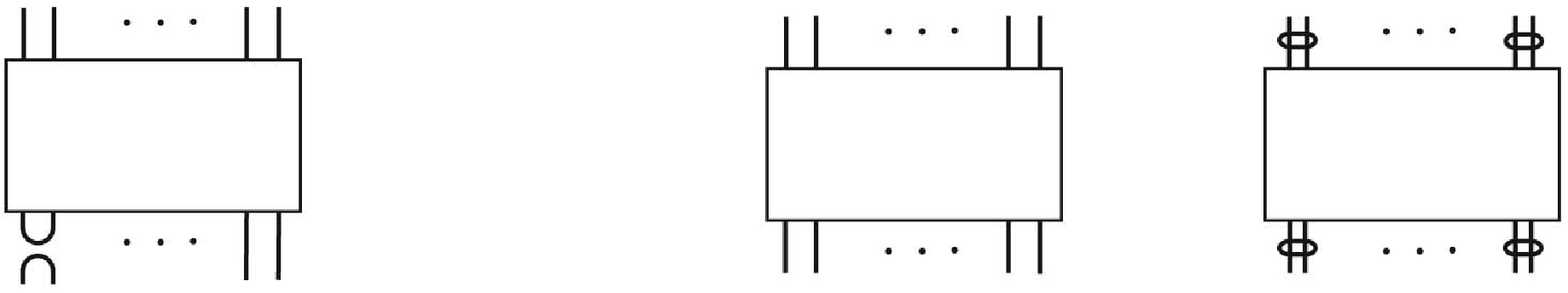}
    \put(-285,26){$P^{(2m)}$}
    \put(-138,25){$P^{(2m)}$}
    \put(-42,25){$P^{(2m)}$}
    \put(-233,27){$=\;0,$}
    \put(-83,27){$=$}
     \caption{}
     \label{fig:JW2}
\end{figure}

{\em Proof}. Any element ${\mathcal E}$ of $TL_n$ obeying
$E_j{\mathcal E}={\mathcal E}E_j=0$ for all $j<n$ is said to be ``killed
by turnbacks''. It follows that any such element
obeys $P^{(2)}_j{\mathcal E}P^{(2)}_j={\mathcal E}$.
As is straightforward to check using the recursion relation, the
Jones-Wenzl projectors are killed by turn-backs, as illustrated
in the left of figure \ref{fig:JW2}.
Thus
$$P^{(2m)}= P^{(2)}_1P^{(2)}_3\dots P^{(2)}_{2m-1} P^{(2m)}
P^{(2)}_1P^{(2)}_3\ \dots P^{(2)}_{2m-1}\
 .$$ We can thus pair up the
strands above and below the Jones-Wenzl projector, and
replace each pair with $P^{(2)}_j$, as illustrated on the right of
figure \ref{fig:JW2}.

\begin{figure}[ht]
\includegraphics[height=1.7cm]{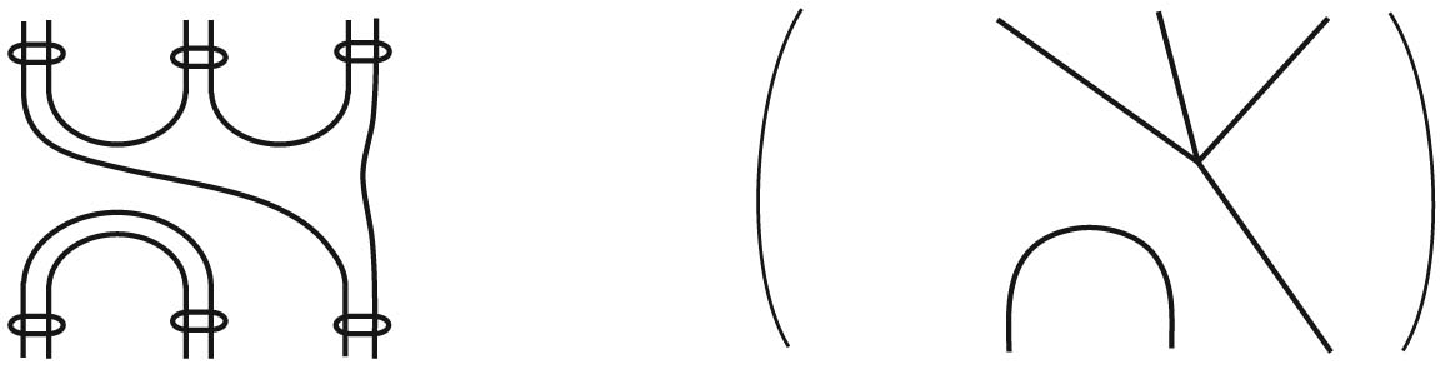}
{\small     \put(-125,22){$=\; {\Phi}$}
    \put(-88,22){$d^{-1}\,\cdot$}}
     \caption{}
\label{fig:image}
\end{figure}

Consider the projector $P^{(2m)}$ as a linear combination of curve
diagrams in the rectangle.  Each pair of external strands of $P^{(2m)}$
is replaced with $P^{(2)}$, so all external strands in the Jones-Wenzl
projector in $TL_{2m}$ correspond to lines in the chromatic
algebra. Since each diagram (additive generator of the Temperley-Lieb algebra)
consists of disjoint embedded curves, it follows from the definition of ${\Phi}$
(figure \ref{fig:anotherphi}) that each individual term
in the expansion of $P^{(2m)}$, with external strands paired up and
replaced with $P^{(2)}$, is in the image of ${\Phi}$, see figure \ref{fig:image} for an
example.
Thus the linear
combination, $P^{(2m)}\in TL_{2m}$, is indeed in the image of the
chromatic algebra ${\mathcal C}_m$. \qed

\begin{figure}[ht]
\includegraphics[width=11.5cm]{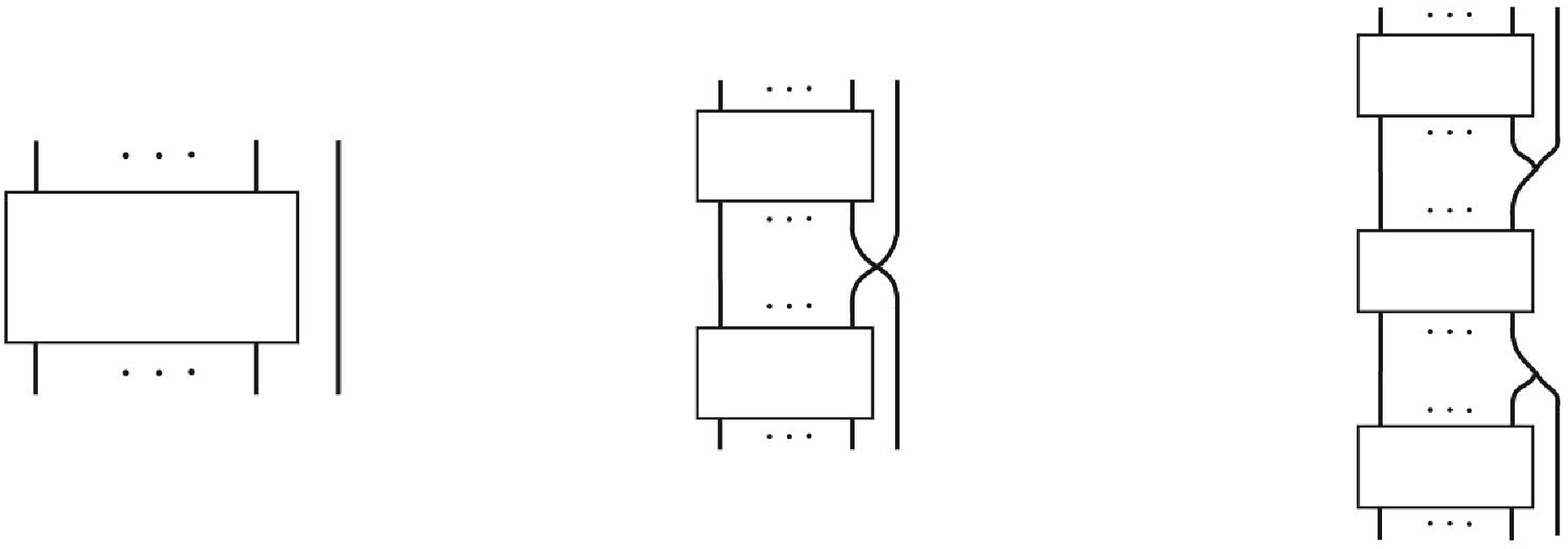}
    \put(-377,55){$\overline P^{(2m)}\, =$}
    \put(-314,55){$\overline P^{(2m-2)}$}
    \put(-246,55){$-\, \frac{1}{d} \cdot \frac{\Delta_{2m-3}}{\Delta_{2m-2}}$}
{\scriptsize
    \put(-180,32){$\overline P^{(2m-2)}$}
    \put(-180,78){$\overline P^{(2m-2)}$}}
    \put(-133,55){$-\; \frac{1}{d}\cdot \frac{({\Delta}_{2m-3})^2}{{\Delta}_{2m-1}{\Delta}_{2m-2}}$}
{\scriptsize
    \put(-41,95){$\overline P^{(2m-2)}$}
    \put(-41,53){$\overline P^{(2m-2)}$}
    \put(-41,13){$\overline P^{(2m-2)}$}}
     \caption{A recursive formula for the pull-back $\overline P^{(2m)}$ of the Jones-Wenzl projector
     $P^{(2m)}$ in the chromatic algebra.}
\label{fig:recursive}
\end{figure}

Applying the recursive formula for the Jones-Wenzl projector (figure \ref{fig:JonesWenzl}) twice and
using the technique in the proof of lemma \ref{imagelemma}, one gets the formula for
the preimage of $P^{(2m)}$ in the chromatic algebra:

\begin{lem} \label{recursive} \sl
The pull-back $\overline P^{(2m)}$ of the Jones-Wenzl projector $P^{(2m)}$ to ${\mathcal C}_m$ is given by the
recursive formula in figure \ref{fig:recursive}.
\end{lem}

The base of this recursion is $\overline P^{(2)}$ which is now just a
single strand.  As mentioned above, the commutativity of diagram
(\ref{another commutative}) implies that the pullback of the trace
radical in $TL^{\sqrt{Q}}_{2m}$ is in the trace radical of ${\mathcal
  C}_m^{Q}$.  Therefore lemma \ref{recursive} establishes
a chromatic polynomial identity generalizing
(\ref{Tutte's identity3}) for each
value of $Q$ obeying
$$Q=4\cos^2\left(\frac{\pi j}{2m+1}\right) = 2+ 2\cos\left(\frac{2\pi
  j}{2m+1}\right)$$ with $j<2m$. This generalized identity may be
generated explicitly by using the recursion relation for $\overline
P^{(2m)}$ in figure \ref{fig:recursive}.  Specifically, the relation
$\overline P^{(2m)}=0$ preserves the chromatic polynomial of the dual
graphs, and replacing each graph in the relation $\overline
P^{(2m)}=0$ by its dual gives a generalization of Tutte's relation
(\ref{Tutte's identity3}).  For example, using the recursive formula
shows that (\ref{Tutte's identity4}) (checked directly in section
\ref{sec:TL} by showing ${\Phi}$ maps it to $P^{(4)}$) is equivalent in
the chromatic algebra to setting $\overline P^{(4)}=0$.

\subsection{Other values of Q}
We have used the Jones-Wenzl projectors labeled by even integers to find chromatic
identities for half the values of $Q$ in (\ref{generalizedBeraha}). In
this subsection, we show how to use the projectors $P^{(2m-1)}$ to
find chromatic identities valid at $Q$ in (\ref{generalizedBeraha})
for odd $n$.

We start by indicating a direct generalization of the argument above, starting
with the Jones-Wenzl projector $P^{(2m-1)}\in TL_{2m-1}$ and getting a relation
in ${\mathcal C}_m$ at the corresponding value of $Q$. The drawback of this approach is that
it involves a choice of including $P^{(2m-1)}$ into $TL_{2m}$. Further below we pass from algebraic
to categorical setting to get a unique chromatic relation corresponding to $P^{(2m-1)}$.

\begin{figure}[ht]
\includegraphics[height=2cm]{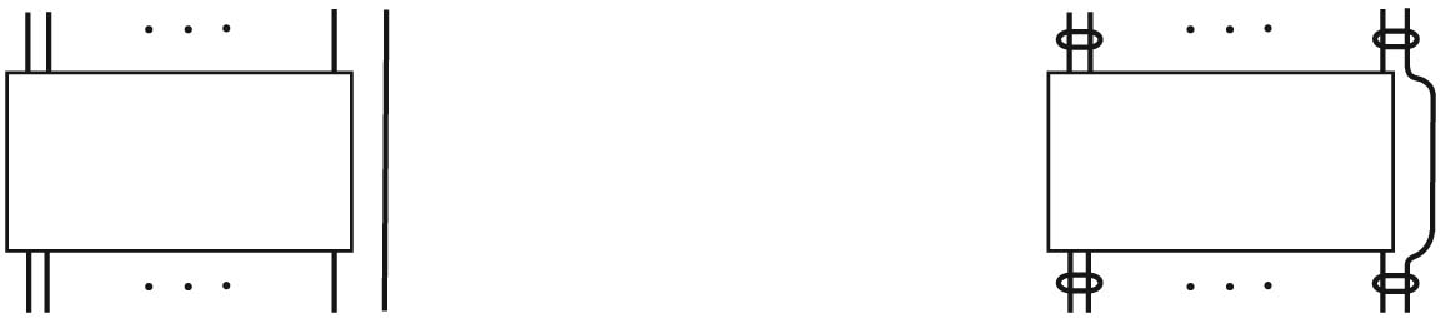}
    \put(-315,23){$P_1^{(2m-1)}\; =$}
    \put(-240,23){$P^{(2m-1)}$}
    \put(-131,23){$R^{(2m-1)}\; =$}
    \put(-57,23){$P^{(2m-1)}$}
     \caption{}
     \label{fig:JW7}
\end{figure}

To find chromatic identities for
even $n$, we paired up the lines in $P^{(2m)}\in TL_{2m}$
and showed that this projector could be pulled back to the
chromatic algebra. Since $P^{(2m-1)}\in TL_{2m-1}$ acts on odd number
of lines, finding the generalization of lemma \ref{imagelemma}
requires a little more work. The simplest way is to map
$P^{(2m-1)}$ to an element in $TL_{2m}$
by adding a non-intersecting strand at the
right; we label this as $P_1^{(2m-1)}\in TL_{2m}$, figure \ref{fig:JW7}.
$P_1^{(2m-1)}$ is killed by turnbacks $E_j$ for $j<2m-1$, but to be
able to pair up all $2m$ lines, we need to define
the element $R^{(2m-1)}\in TL_{2n}$ as
\begin{equation}
R^{(2m-1)}\equiv  P^{(2)}_{2m-1} P^{(2m-1)}_1 P^{(2)}_{2m-1}=
P^{(2)}_1P^{(2)}_3\dots P^{(2)}_{2m-1} P^{(2m-1)}_1
P^{(2)}_1P^{(2)}_3\ \dots P^{(2)}_{2m-1}
\ ,
\label{Rdef}
\end{equation}
figure \ref{fig:JW7}. $R$ is no longer a projector, but is killed by all turnbacks with
$j\le 2m-1$. Therefore we can pair up the strands as before, and
a rerun of the argument proving the previous lemma shows that:
\begin{lem} \sl The element $R^{(2m-1)}\in TL_{2m}$ is in the image of the chromatic algebra: $R^{(2m-1)}\in {\Phi}({\mathcal C}_m)$,
so its pullback at the corresponding value of $Q$ is in the trace
radical of the chromatic algebra.
\end{lem}
$R^{(2m-1)}$ is in the ideal generated by $P^{(2m-1)}$,
thus it is in the trace radical when  $Q=d^2$ obeys (\ref{generalizedBeraha})
for odd $n=2m-1$. A chromatic identity then follows by
taking the sum of the dual graphs of the pictures for the pull-back ${\Phi}^{-1}(R^{(2m-1)})$.

This construction of an ideal in ${\mathcal C}_m$ at
$\sqrt{Q}=2\cos(\pi/(2m))$ is not unique. Other elements of ${\mathcal C}_m$ at this value of $Q$
can be set to zero as well, although it is not clear to us if any of
these result in new chromatic identities, or simply rotations of each
other.

In order to pull back $P^{(2m-1)}$ to ${\mathcal C}_{m}$, we needed to
add an extra line so that we could form $m$ pairs on the bottom and
$m$ on the top. To pull back $P^{(2m-1)}$ itself, and thus to
find a unique chromatic identity, we do not add an
extra line, but instead pair up a line from the top of $P^{(2m-1)}$
with one from the bottom. The pullback of this object no longer lives
in the chromatic algebra, but rather the {\em chromatic category},
briefly discussed at the end of subsection 4.2.  The
three equivalence relations (including the contraction-deletion rule) $(1)-(3)$
at the beginning of subsection 4.1 serve as defining relations in this
category as well. Finding this pull-back
results in a chromatic identity involving the duals of graphs on a disk
with $2m-1$ external strands in total. To be more precise, we will consider
graphs in a disk $D$ with $n$ fixed points on the boundary of $D$, modulo the
relations $(1)-(3)$ in section 4.2. Moreover, the marked points on the boundary are
numbered $1$ through $n$.
In the categorical language, these are morphisms: elements of $Hom(0,n)$. Given two such
graphs $a, b$, their inner product is computed by reflecting $b$ and gluing the two disks
so the numbered points on the boundary are matched. Then $\langle a,b\rangle$ is given
by the evaluation of the resulting graph in the sphere (equal to $Q^{-1}$ times the
chromatic polynomial of the dual graph.) In this setting, the trace radical is replaced
by the ideal of {\em negligible morphisms}, see \cite{F} for
more details.
We extend the map $\Phi$
to the chromatic category (with values in the Temperley-Lieb category) in the obvious way:
it is defined by replacing each line and vertex
by the doubled lines in $TL$, as
illustrated in figures \ref{fig:phi1} and \ref{fig:anotherphi}.

\begin{figure}[ht]
\includegraphics[height=2cm]{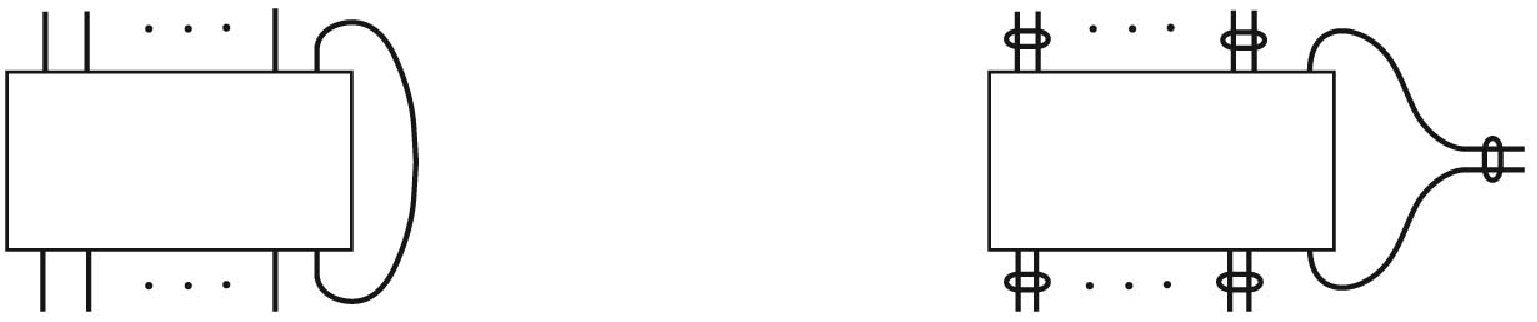}
    \put(-327,23){$tr_n\, P^{(n)}\, =$}
    \put(-251,23){$P^{(n)}$}
    \put(-155,23){$P^{(2m+1)}=$}
    \put(-85,23){$P^{(2m+1)}$}
     \caption{The identity on the right holds in the Temperley-Lieb category at $d=2\, \cos({\pi}j/(2m+2)).$}
     \label{fig:JW8}
\end{figure}

We start by showing when the Jones-Wenzl projectors are killed by
``end turnbacks''. An end turnback on the right in $TL_n$ is the
partial trace $\hbox{tr}_n: TL_n\to TL_{n-1}$, defined by connecting
just the $n$th strand on the bottom to the $n$th strand on the top,
figure \ref{fig:JW8} (left). Likewise, an end
turnback on the left is the partial trace $\hbox{tr}_1$. From the
recursion relation pictured in (\ref{fig:JonesWenzl}), we find
$$\hbox{tr}_n(P^{(n)})= \left(d -
\frac{\Delta_{n-2}}{\Delta_{n-1}}\right) P^{(n-1)} =
\frac{\Delta_n}{\Delta_{n-1}} P^{(n-1)}
$$
When $d$ takes on the special values (\ref{special d}), $\Delta_n=0$,
so only at these values is
$P^{(n)}$ killed by the end turnback.

When $P^{(n)}$ is killed by all turnbacks including the end one,
it can be pulled back to the chromatic category for any $n$, odd or
even. The idea is the same as for the algebras. We pair up the lines
in the Temperley-Lieb category, and then replace each pair with a
single line in the chromatic category. Figure \ref{fig:JW8} (right) illustrates the
odd case $n=2m+1$. Note that while lemma \ref{imagelemma} is true for any value of $Q$,
this categorical analogue for $n=2m+1$
holds only at the special value of $Q=2+2\cos(2{\pi}j/(2m+2))$ (and the corresponding
value of the parameter $d=\sqrt Q$)!

\begin{figure}[ht]
\includegraphics[height=3.2cm]{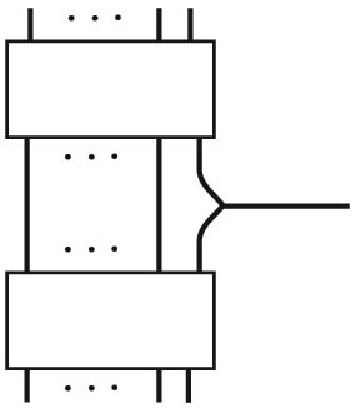}
    \put(-183,42){$\overline P^{(2m+1)}\; =\; -\; \frac{{\Delta}_{2m-1}}{{\Delta}_{2m}}$}
    \put(-66,65){$\overline P^{(2m)}$}
    \put(-66,15){$\overline P^{(2m)}$}
     \caption{A recursive formula for the pull-back $\overline P^{(2m+1)}$ of the Jones-Wenzl projector
     $P^{(2m+1)}$ in the chromatic category.}
\label{fig:recursive2}
\end{figure}

Considering the recursive formula (figure \ref{fig:JonesWenzl}) for the Jones-Wenzl projector $P^{(2m+1)}$, one gets
the following formula for its pullback:

\begin{lem} \label{recursive2} \sl
The pull-back $\overline P^{(2m+1)}$ of the Jones-Wenzl projector $P^{(2m+1)}$ to the chromatic category is given by the
recursive formula in figure \ref{fig:recursive2}.
\end{lem}

Together with lemma \ref{recursive}, this gives a recursive formula
for the pullback $\overline P^{(n)}$ of the Jones-Wenzl projector
$P^{(n)}$ for all values of $n$.  Considering the dual graphs for the
graphs in the relation $\overline P^{(n)}=0$, one gets a chromatic
polynomial relation for each value of $Q$ in
(\ref{generalizedBeraha}).

To give an example, $\overline{P}^{(3)}$ simply is a trivalent
vertex. This indeed is in the trace radical when $Q=2$, as is easy
to see by reverting to the original definition of the chromatic
polynomial: the dual graph of a trivalent vertex
is a triangle, and any graph containing a triangle cannot be colored with two
colors.

\medskip

{\bf Acknowledgments.} We would like to thank Mike Freedman and Kevin
Walker for discussions on the Temperley-Lieb category and related
questions in quantum topology, as well as for introducing us to the
golden identity.

\end{document}